\documentclass[12pt]{article}

\usepackage{amsfonts,amssymb,amstext,graphics,amsmath,color}

\usepackage{amsmath,amssymb,amsbsy,upref}

\usepackage[english]{babel}
\usepackage{graphicx}
\usepackage{graphics} 
\usepackage{epsfig} 
\usepackage{color}
\usepackage{amsthm}  

\numberwithin{equation}{section}

\newtheorem{theorem}{Theorem}[section]
\newtheorem{lemma}[theorem]{Lemma}
\newtheorem{proposition}[theorem]{Proposition}
\newtheorem{corollary}[theorem]{Corollary}

\newtheorem{remark}[theorem]{Remark}

\usepackage{titlesec}
\titleformat{\section}
{\bf\large}{\thesection}{7pt}{\bf\large}
\titleformat{\subsection}
{\bf}{\thesubsection}{5pt}{\bf}
\titleformat{\subsubsection}
{\bf}{\thesubsubsection}{3pt}{\bf}

\newcommand{\ee}{\varepsilon}

\renewcommand{\k}{\kappa}

\newcommand{\lie}{\mathop{\displaystyle {\stackrel{\scriptstyle{\varepsilon \to 0}}{\hbox to 40 pt {\rightarrowfill}}  } }}

\newcommand{\lil}{\mathop{\displaystyle {\stackrel{\scriptstyle{l\to \infty}}{\hbox to 40 pt {\rightarrowfill}}  } }}

\begin{document}

\title{\bf\Large Asymptotic stability of the spectrum of a parametric \\ family of homogenization problems associated \\ with a perforated waveguide\footnote{This is the pre-peer reviewed version of the following article: [Delfina Gómez, Sergei A. Nazarov, Rafael Orive-Illera, María-Eugenia Pérez-Martínez. Asymptotic stability of the spectrum of a parametric family of homogenization problems associated with a perforated waveguide, Mathematische Nachrichten, 2023, DOI: 10.1002/mana.202100589], which will be published in final form at [https://doi.org/10.1002/mana.202100589]. This article may be used for non-commercial purposes in accordance with Wiley Terms and Conditions for Use of Self-Archived Versions.}}
\author{\normalsize{Delfina G{\'o}mez$^a$, Sergei A. Nazarov$^b$, Rafael Orive-Illera$^{c,d}$, Maria-Eugenia P\'erez-Mart\'inez$^e$ } \vspace{0.2cm}\\
\small{$^a$ Departamento de Matem\'aticas, Estad\'{\i}stica y Computaci\'on, Universidad de Cantabria,} \vspace{-0.1cm} \\
\small{Santander, Spain, \texttt{gomezdel@unican.es}} \vspace{0.1cm} \\
\small{$^b$ Institute for Problems in Mechanical Engineering of the Russian Academy of Sciences,} \vspace{-0.1cm}\\
\small{Saint Petersburg,  Russia, \texttt{srgnazarov@yahoo.co.uk}}\vspace{0.1cm}\\
\small{$^c$ Instituto de Ciencias Matem\'aticas, CSIC-UAM-UC3M-UCM, Madrid, Spain,}\vspace{-0.1cm}\\
\small{$^d$ Departamento de Matem\'aticas, Universidad Aut\'onoma de Madrid,}\vspace{-0.1cm}\\
\small{Madrid, Spain, \texttt{rafael.orive@icmat.es}}\vspace{0.1cm}\\
\small{$^e$ Departamento de Matem\'atica Aplicada y Ciencias de la Computaci\'on, Universidad de Cantabria,}\vspace{-0.1cm}\\
\small{Santander, Spain, \texttt{meperez@unican.es}}
}

\date{}

\maketitle

\vspace*{-0.7cm}
\begin{center}
\begin{minipage}{6in}
\thispagestyle{empty} \setlength{\baselineskip}{6pt}
{\small {\bf Abstract}:
In this paper, we provide uniform bounds for convergence rates  of the   low frequencies of a parametric family of  problems for the Laplace  operator posed on a  rectangular perforated domain of the plane of height $H$.  The perforations are periodically placed along the ordinate  axis at a  distance  $O(\varepsilon)$ between them, where $\varepsilon$ is a parameter that converges towards zero. Another parameter $\eta$, the  Floquet-parameter,     ranges in the interval $[-\pi,\pi]$. The boundary conditions are quasi-periodicity conditions on the lateral sides of the rectangle and Neumann over the rest.
We  obtain  precise bounds for convergence rates which are uniform on  both parameters $\varepsilon$ and $\eta$ and strongly depend on $H$. As a model problem associated with a waveguide, one of the main difficulties in our analysis comes  near the nodes of the limit dispersion curves.

\vspace*{0.2cm}

\vspace*{0.2cm}
}
\end{minipage}
\end{center}

\vspace*{0.5cm}

\setcounter{equation}{0} \setcounter{section}{0}
\renewcommand{\theequation}{\arabic{section}.\arabic{equation}}
\section{Introduction.}\label{sec1}
In this section, we formulate the  spectral perturbation problem which constitutes a homogenization problem for the Laplacian in a periodically perforated  rectangular  domain, $\ee H $ being the period, with $\ee\ll 1$ and $H>0$; $H$  is the height of the rectangle. The so-called {\it Floquet-parameter}  $\eta\in [-\pi,\pi]$ arises on  the  lateral boundary conditions.  On the rest of the boundary, Neumann   conditions are imposed. The perforations are periodically placed along the ordinate axis.
In Section \ref{sec1.2},  we introduce the homogenized problem which is  still a $\eta$-dependent problem. Some background on the problem and the structure of the paper are in Section \ref{sec1.3}. Our aim is to study the asymptotic behavior of the spectrum, as $\ee\to 0$, in its dependence on the other parameters   $\eta$ and $H$.

\medskip

 \subsection{The parametric family of  homogenization spectral problems. }\label{sec1.1}
Let  $H$ be a positive parameter and let  $\varpi^0$  be the  rectangle    \begin{equation}\label{rectangulo}\varpi^0=\{x=(x_1,x_2)\in{\mathbb R}^2:\,|x_1|<1/2,\quad x_2\in(0,H)\}.\end{equation}

Let $\omega$ be a domain in the
plane ${\mathbb R}^2$ which is bounded by a smooth simple closed curve
$\partial \omega$ and has the compact closure
$\overline{\omega}=\omega\cup\partial\omega$ inside $\varpi^0$. We
introduce the perforated domain  $\varpi^\varepsilon$, see Figure~\ref{fig1}, a), obtained from $\varpi^0$  by removing the family of holes
 \begin{equation}\begin{array}{c}
 \omega^\varepsilon(k)=\{x:\,\varepsilon^{-1}(x_1 ,x_2-  \varepsilon k { \color{black} H})\in\omega\},\quad
 k=0,\dots,N-1,
 \end{array}\label{(2)}\end{equation}
which are distributed periodically along the  line   $x_1=0$. Each hole is homothetic to $\omega$ of ratio $\varepsilon$ and translation of $\varepsilon\omega=\omega^\varepsilon(0).$ Namely,
\begin{equation}\begin{array}{c}\displaystyle
\varpi^\varepsilon=\varpi^0\setminus\overline{\omega^\varepsilon}\quad \mbox{\rm
where}\quad\omega^\varepsilon=
\bigcup\limits_{k=0}^{N-1} \omega^\varepsilon(k).
\end{array}\label{(3)}\end{equation}
Here,   $\varepsilon$  is a small positive
parameter and $N$ is a big natural
number, both related by  {$N=\ee^{-1}$.
 The period   is $\varepsilon H $ with $\varepsilon\ll 1.$

\begin{figure}[t]
\begin{center}
\resizebox{!}{2.5cm}{\includegraphics{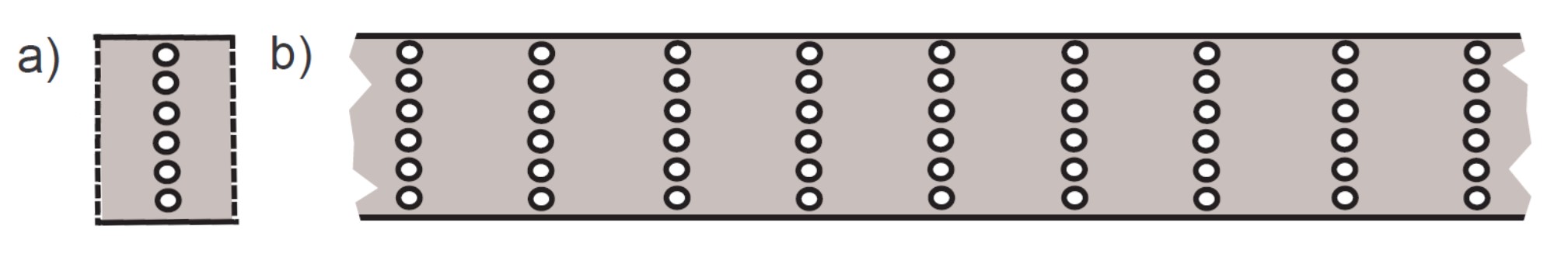}}
\caption{{\bf a)} The perforated domain  $\varpi^\ee$.  {\bf b)} The perforated strip $\Pi^\ee$ .}
\label{fig1}
\end{center}
\end{figure}

In the domain $\varpi^\ee$, we   consider the spectral problem  consisting in   the partial differential equation
 \begin{equation}\begin{array}{c}
-\Delta{\color{black}_x} U^\varepsilon(x;\eta)=\Lambda^\varepsilon(\eta)U^\varepsilon(x;\eta),
\,\, x\in\varpi^\varepsilon,
\end{array}\label{(11)}\end{equation}
the quasi-periodicity conditions on the vertical sides of $\varpi^\ee$
 \begin{equation}\begin{array}{c}\displaystyle
U^\varepsilon\Big(\frac{1}{2},x_2;\eta\Big)=e^{{\rm i}\eta}U^\varepsilon\Big(-\frac{1}{2},x_2;\eta\Big),
\,\, x_2\in(0,H),
\\\\\displaystyle
\frac{\partial U^\varepsilon}{\partial x_1}\Big(\frac{1}{2},x_2;\eta\Big)=e^{{\rm i}\eta}
\frac{\partial U^\varepsilon}{\partial x_1}\Big(-\frac{1}{2},x_2;\eta\Big),
\,\, x_2\in(0,H),
\end{array}\label{(12-13)}\end{equation}
and the Neumann condition on the remaining part
of the boundary of $\varpi^\varepsilon$
 \begin{equation}\begin{array}{c}
\partial_\nu U^\varepsilon(x{\color{black};\eta})=0,\,\, x\in \{x\in\partial
\varpi^\varepsilon:\,|x_1|<1/2\}.
\end{array}\label{(14)}\end{equation}
Here,  $\eta\in[-\pi,\pi]$ is the Floquet-parameter,
$\partial_\nu$ denotes the directional derivative along the outward
normal $\nu$.
 $\Lambda^\varepsilon(\eta)$ and
$U^\varepsilon(\cdot;\eta)$, respectively, denote
the eigenvalues and eigenfunctions which depend on both the perturbation parameter and the Floquet-parameter. We address the  asymptotic behavior of $(\Lambda^\varepsilon(\eta),\,  U^\varepsilon(\cdot;\eta))$ as $\ee\to 0$.

The variational formulation of the problem \eqref{(11)}--\eqref{(14)} reads (see, e.g., \cite{Lad}): Find  $\Lambda^\ee(\eta)$ and $ U^\varepsilon(\cdot;\eta)\in H^{1,\eta}_{per}(\varpi^\varepsilon)$, $U^\varepsilon(\cdot;\eta)\not=0 $, satisfying
\begin{equation}\begin{array}{c}
 (\nabla_x U^\varepsilon(\cdot;\eta),\nabla_x V^\varepsilon)_{\varpi^\varepsilon}=\Lambda^\varepsilon(\eta)
 (U^\varepsilon(\cdot;\eta),V^\varepsilon)_{\varpi^\varepsilon}\quad \forall V^\varepsilon\in H^{1,\eta}_{per}(\varpi^\varepsilon),
\end{array}\label{(16)}\end{equation}
where $H^{1,\eta}_{per}(\varpi^\varepsilon)$ is the Sobolev space of
functions  in $H^{1}(\varpi^\varepsilon)$ satisfying the quasi-periodi\-city conditions \eqref{(12-13)}, while $(\cdot,\cdot)_{\varpi^\varepsilon}$ stands for the
scalar product in  $L^2(\varpi^\varepsilon)$.

  In view of the compact embedding
$H^1(\varpi^\varepsilon) \subset L^2(\varpi^\varepsilon)$
problem \eqref{(16)}  has a discrete spectrum constituting the
unbounded monotone sequence of eigenvalues (cf.  \cite[\S\,10.2]{BiSo} and \cite[\S\,4.5]{SHSP}),
 \begin{equation}\begin{array}{c}
0\leq\Lambda^\varepsilon_1(\eta)\leq\Lambda^\varepsilon_2(\eta)\leq\dots\leq
\Lambda^\varepsilon_p(\eta)\leq\dots\to+\infty, \quad {\color{black} \mbox{ as } p\to +\infty},
\end{array}\label{(17)}\end{equation}
which  are  repeated  according to their multiplicities. Also the corresponding eigenfunctions   $\{{ U_p^\ee(\cdot;\eta)} \}_{p=1}^\infty$ are chosen to form an orthonormal basis in $L^2(\varpi^\ee)$.

Furthermore, the functions
 \begin{equation}\begin{array}{c}
[-\pi,\pi]\ni \eta\,\,\mapsto\,\, \Lambda^\varepsilon_p(\eta),\quad
p\in{\mathbb N},
\end{array}\label{(18)}\end{equation}
are continuous and $2\pi$-periodic. This last assertion is due to the fact that   problem   \eqref{(11)}-\eqref{(14)} is the model problem associated with a waveguide: a periodically perforate Neumann strip recently considered in the literature (cf. \eqref{(4-5)}, Figure \ref{fig1} b) and \cite{nuestroNeumann}).  For the sake of completeness, in order to   outline the interest of the problem under consideration \eqref{(11)}-\eqref{(14)}, as well as its properties  we  introduce briefly  this waveguide in Section \ref{sec1.3}.

\medskip

\subsection{ The parametric family of  homogenized problems. }\label{sec1.2}
By analogy with the homogenization of perforated domains along manifolds with   Neumann boundary conditions (see, e.g., \cite{nuestroOleinik}), we easily see that the homogenized problem is a spectral problem for the Laplacian, posed in the rectangle  $\varpi^0$, cf. \eqref{rectangulo}, with the stable Neumann \eqref{(30)} (quasi-periodicity \eqref{(25)}, respectively) conditions on the horizontal  (vertical, respectively) sides of the rectangle. For the readers convenience, we introduce here this problem and provide the  explicit formulas for the eigenvalues and eigenfunctions in terms of  the parameters $\eta$ and $H$.  The    convergence of the spectrum of \eqref{(11)}-\eqref{(14)}, as $\ee\to 0$, will be outlined in Section \ref{sec2} (cf. Corollary \ref{Lemma_convergence}) as a consequence of a more general result (cf. Theorem \ref{corollary_convergence}).

In $\varpi^0$,  we set   the {\it homogenized spectral problem}
 \begin{equation}\begin{array}{c}
-\Delta{\color{black}{_x}} U^0(x;\eta)=\Lambda^0(\eta)U^0(x;\eta), \,\, x\in\varpi^0,
\end{array}\label{(29)}\end{equation}
 \begin{equation}\begin{array}{c}\displaystyle
\frac{\partial U^0}{\partial x_2}(x_1,0{\color{black};\eta})=\frac{\partial
U^0}{\partial x_2}(x_1,H{\color{black};\eta})=0,\,\,
x_1\in\Big(-\frac{1}{2},\frac{1}{2}\Big),
\end{array}\label{(30)}\end{equation}
\begin{equation}\begin{array}{c}\displaystyle
U^0\Big(\frac{1}{2},x_2;\eta\Big)=
e^{{\rm i}\eta}U^0\Big(-\frac{1}{2},x_2;\eta\Big),
\\\\\displaystyle
 \frac{\partial U^0}{\partial
x_1}\Big(\frac{1}{2},x_2;\eta\Big)= e^{{\rm i}\eta}\frac{\partial
U^0}{\partial x_1}\Big(-\frac{1}{2},x_2;\eta\Big), \,\, x_2\in(0,H),
\end{array}\label{(25)}\end{equation}
 where $\Lambda^0(\eta)$ and $U^0(\cdot;\eta)$  denote  the eigenvalues  and corresponding eigenfunctions.

\begin{figure}[t]
\begin{center}
\resizebox{!}{8cm} {\includegraphics{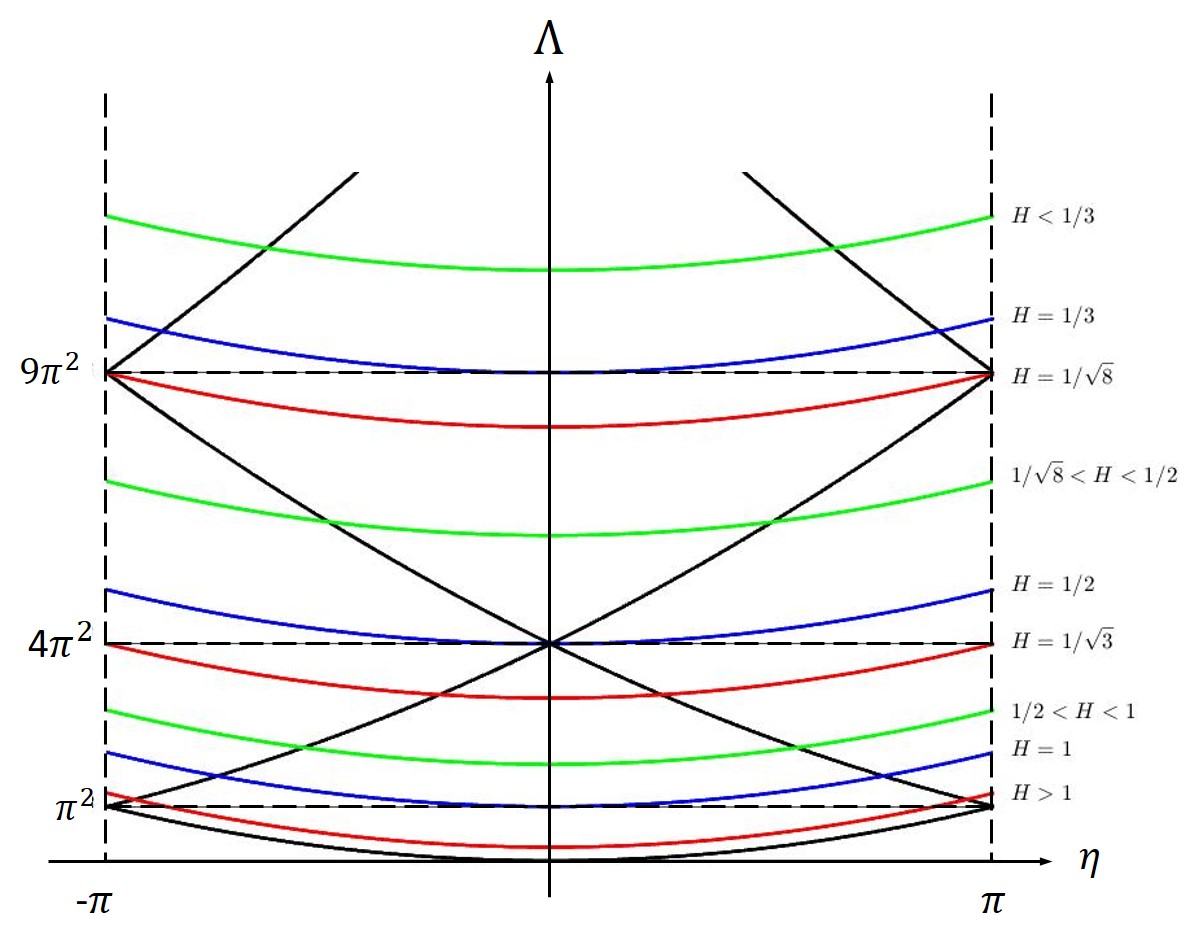}}
\caption{The first limit dispersion curves  for different values of $H$. The first one is $\Lambda=\eta^2$. The parabolic curves translations of this one along the $\Lambda$-axis are the curves  $\Lambda =\eta^2+\pi^2  H^{-2}$ depending on $H$. The other pieces of parabola are the translations along the $\eta$-axis, namely, $\Lambda=(\eta\pm 2\pi j)^2$, while $j=1,2$.}
\label{fig_variosH}
\end{center}
\end{figure}

 Problem \eqref{(29)}-\eqref{(25)}  admits the variational formulation: Find  $\Lambda^0(\eta)$ and  $U^0(\cdot;\eta)\not=0 $, $ U^0(\cdot;\eta)\in H^{1,\eta}_{per}(\varpi^0)$,   satisfying
\begin{equation}
\begin{array}{c}
  (\nabla_x {  U}^0(\cdot;\eta) ,\nabla_x  V)_{\varpi^0}=  \Lambda^0 (\eta)({
U}^0(\cdot;\eta),V)_{\varpi^0}
  \quad \forall   V\in H^{1,\eta}_{per}(\varpi^0).
\end{array}\label{(16limit)}\end{equation}
As is well known,
 it  has a discrete spectrum, which we can write in an increasing order
 \begin{equation}\begin{array}{c}
0\leq\Lambda^0_1(\eta)\leq\Lambda^0_2(\eta)\leq\dots\leq
\Lambda^0_p(\eta)\leq\dots\to+\infty,\quad {\color{black} \mbox{ as } p\to +\infty},
\end{array}\label{(S2)}\end{equation}
where the convention of repeated eigenvalues has been adopted. Also the corresponding eigenfunctions   $\{{ U^0_p(\cdot;\eta)}\}_{p=1}^\infty$ are chosen to form an orthonormal basis in $L^2(\varpi^0)$.
Furthermore, the functions
 \begin{equation}\begin{array}{c}
[-\pi,\pi]\ni \eta\,\,\mapsto\,\, \Lambda^0_p(\eta),\quad
p\in{\mathbb N},
\end{array}\label{(18m)}\end{equation}
are continuous and $2\pi$-periodic.

Explicit formulas for eigenvalues and eigenfunctions of  \eqref{(29)}-\eqref{(25)}  are obtained by separation of variables:
\begin{equation}\label{(31)}
\Lambda^0_{jk}(\eta)=(\eta+2\pi j)^2+\pi^2\frac{k^2}{H^2} \, ,
\qquad  U^0_{jk}(x;\eta)=e^{{\rm i}(\eta+2\pi j)x_1}\cos\Big(\pi
k\frac{x_2}{H}\Big),
\quad
 j\in{\mathbb Z},\,\,k\in{\mathbb N}_0={\mathbb N}\cup\{0\}.
\end{equation}
It should be mentioned that re-numeration of the eigenvalues in
\eqref{(31)} is needed to compose the monotone sequence
\eqref{(S2)}.

Note that formulas \eqref{(31)} are of great interest to draw the {\it limit dispersion curves}  for different values of $H$. Recall that these curves are the graphs of $\Lambda_p^0(\eta)$, for $\eta\in [-\pi,\pi]$. Figure \ref{fig_variosH} shows the great variety of possibilities of behaviors of such   curves depending on the value of $H$. Also, along with Figures \ref{fig_new} and \ref{fig2},  it gives an idea of what we can expect for the behavior of the perturbed dispersion curves, see, for instance, the second row in  Figure \ref{fig2}.

\begin{figure}[t]
\begin{center}
\resizebox{!}{7cm} {\includegraphics{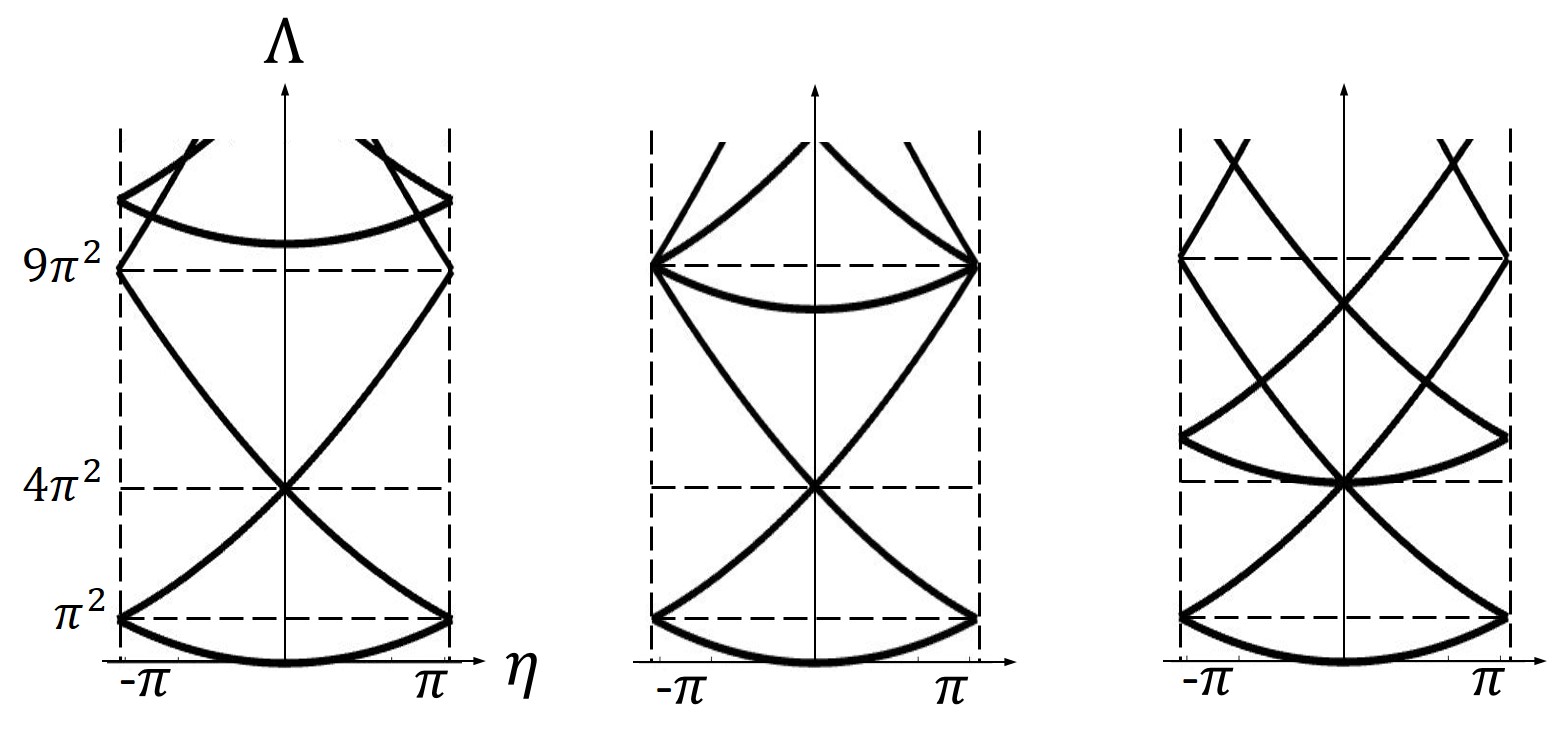}}
\caption{ The  first limit dispersion curves   in the cases  $H<1/3$,  $H=1/\sqrt{8}$ and $H=1/2.$}
\label{fig_new}
\end{center}
\end{figure}

\subsection{The state-of-the-art and the structure of the paper}\label{sec1.3}
First, let us introduce  a problem closely related to \eqref{(11)}-{\eqref{(14)}: a Neumann problem for the Laplace operator in a periodically perforate strip with a double periodicity.

\begin{figure}[t]
\begin{center}
\resizebox{!}{4.5cm} {\includegraphics{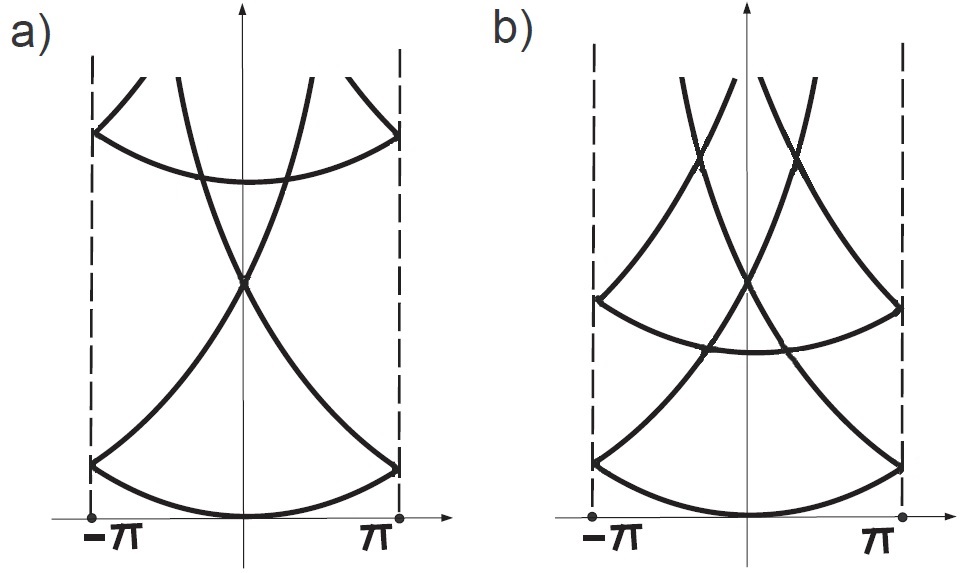}}

\resizebox{!}{4.5cm} {\includegraphics{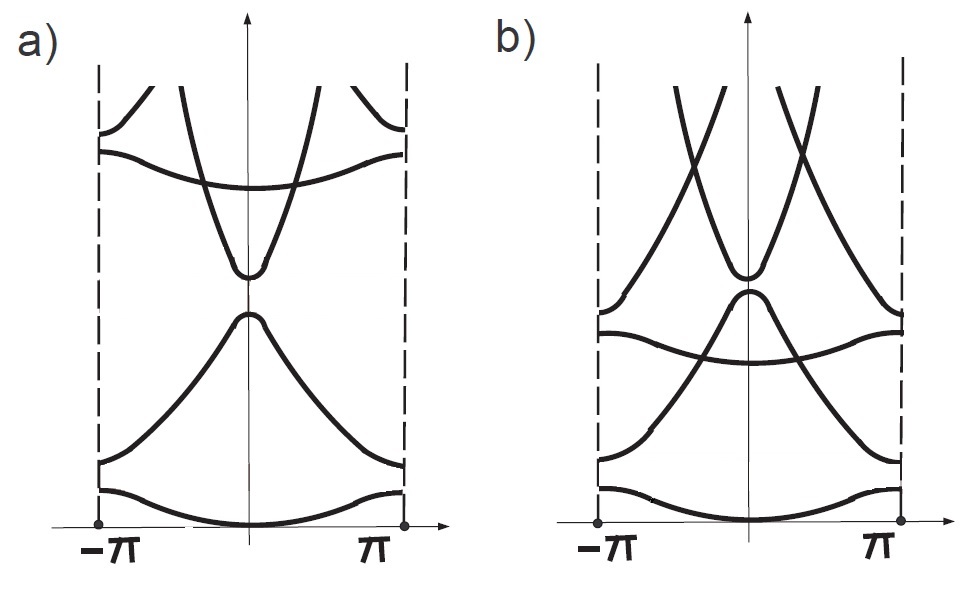}}
\caption{Above,  the  first  six dispersion curves of  the homogenized problem for some values of $H,\,\,$
 a) $H\in(1/\sqrt{8},1/2),\,\,$ b) $H\in(1/2,1)$. Below,  some  possible  dispersion curves  of the perturbed problem for the same values of $H$. }
\label{fig2}
\end{center}
\end{figure}

Extending $\varpi^\ee$ (cf. \eqref{(3)} and Figure~\ref{fig1}}, a)) by periodicity along the $x_1$ axis, we create the unbounded perforated  strip $\Pi^\ee$ (see Figure~\ref{fig1}, b)):
$$
\Pi^\varepsilon=\mathbb{R}\times (0,H)\setminus \bigcup\limits_{j\in{\mathbb Z}}
\bigcup\limits_{k=0}^{N-1} \overline{\omega^\varepsilon(j,k)}$$
where $\omega^\varepsilon(j,k)=\{x:\,\varepsilon^{-1}(x_1-j,x_2-  \varepsilon k  H)\in\omega\} $ with $j\in{\mathbb Z},\,\,k=0,1,\dots,N-1.$
In the
waveguide $\Pi^\varepsilon$, we consider the Neumann spectral  problem
 \begin{equation}\begin{array}{c}
 -\Delta_x u^\varepsilon(x)=\lambda^\varepsilon u^\varepsilon(x),\,\, x\in\Pi^\varepsilon, \vspace{0.1cm} \\ \displaystyle
 \partial_\nu u^\varepsilon(x)=0,\,\, x\in\partial\Pi^\varepsilon.
 \end{array}\label{(4-5)}\end{equation}
Then, applying the Floquet--Bloch--Gelfand  transform  (see, for instance,
\cite{Gel,ReedSimon,Skrig,Kuchbook,NaPl})
 \begin{equation}\nonumber\begin{array}{c}\displaystyle
u^\varepsilon(x)\,\,\mapsto\,\,U^\varepsilon(x;\eta)=
\frac{1}{\sqrt{2\pi}}\,\sum\limits_{p\in{\mathbb Z}}
e^{-{\rm i}p\eta}u^\varepsilon(x_1+p,x_2) ,
\end{array}\label{(9)}\end{equation}
  problem \eqref{(4-5)}  converts into the boundary value problem \eqref{(11)}-\eqref{(14)} in $\varpi^\ee$.

The spectrum of the operator  on the Hilbert space $L^2(\Pi^\varepsilon)$ associated with problem \eqref{(4-5)} is given by
 \begin{equation}\begin{array}{c}
 \sigma^\ee       =\bigcup\limits_{p\in{\mathbb
 N}}\beta_p^\varepsilon\,
 \end{array}\label{(8)}\end{equation}
 where \begin{equation}\begin{array}{c}
 \beta^\varepsilon_p=\{\Lambda^\varepsilon_p(\eta):\,
 \eta\in[-\pi,\pi]\}\subset \overline{{\mathbb R}_+}.
 \end{array}\label{(19)}\end{equation}
As a consequence of  the previously mentioned continuity of $\Lambda^\ee(\eta) $ in  \eqref{(18)},
the sets $\beta^\varepsilon_p$ are closed, connected, and finite segments. Formulas
\eqref{(8)} and \eqref{(19)}  for the spectrum of  problem \eqref{(4-5)} are well-known in the framework of the
Floquet--Bloch--Gelfand theory  (cf. \cite{BLP,ReedSimon,SHT,Skrig,SHSP,Kuchbook,NaPl,CoPlaVa,AllaireConca_JMPA,nuestroNeumann}).

Therefore, to study the asymptotic behavior of the spectrum of \eqref{(11)}--\eqref{(14)} becomes essential to detect the band gap structure of the spectrum    \eqref{(8)}.   Opening  gaps at the low frequency level  have been broached in \cite{nuestroNeumann}  under  certain symmetry   restrictions for the holes (further specifying, the so-called {\it mirror  symmetry}, see \eqref{(symm)} and Figure~\ref{fig-mirror}): roughly speaking, this  involves  controlling the total number of perturbed eigenvalues inside certain boxes of the band   $[-\pi,\pi]\times [0, \infty)$ in the coordinate axis $(\eta,\Lambda)$.  Nevertheless,   the technique in \cite{nuestroNeumann} does not allow
to obtain uniform bounds for discrepancies   between the dispersion curves  $\{\Lambda^\ee_p(\eta): \eta\in[-\pi,\pi]\}$ and $\{\Lambda^0_p(\eta) : \eta\in[-\pi,\pi]\}$ for different values of $p$.
Instead, the technique in the present paper can contribute to discover how to open these gaps.
However, we do not deal with this task but with obtaining the above-mentioned   bounds.

Comparing with former papers in the literature, we  mention \cite{SOME}, \cite{DESOME}, \cite{nuestroNeumann} and \cite{GoNaOrPe_IMSE} as the closest papers. \cite{SOME} and \cite{nuestroNeumann} address opening gaps in a perforated waveguide for the Laplace operator with Dirichlet and Neumann boundary conditions respectively. Asymptotic stability for the spectrum of problem \eqref{(11)}-{\eqref{(14)} but changing the Neumann conditions for the Dirichlet ones has been addressed in \cite{GoNaOrPe_IMSE},  while \eqref{(11)}-{\eqref{(14)} appears as a limit case of the stiff problem considered in  \cite{DESOME}. It happens that for the Dirichlet condition in \cite{SOME} the limit dispersion curves do not depend on $\eta$  and the analysis in \cite{GoNaOrPe_IMSE} becomes much more simplified than for Neumann.

It is worth  mentioning  that  the explicit formulas \eqref{(31)}  allow us to obtain the precise multiplicity of each eigenvalue $\Lambda_p^0(\eta)$ for each $p=1,2,\cdots$ and for each $\eta\in [-\pi,\pi]$. This proves to be essential for our analysis in Sections~\ref{sec3}--\ref{sec4} which  provides the precise number of eigenvalues of the perturbed problem at a distance $O(\ee)$ of the eigenvalues of the homogenized problem (cf.~Figure~\ref{bandas}). It becomes particularly complicated near the points $\eta_{0,p}$ in which the limit  eigenvalues $\Lambda_p^0(\eta)$ change the multiplicity from $1$ to $2$ or $3$ (cf.  Figures~\ref{fig_new}, \ref{inversos} and \ref{bandas} and Remark \ref{rmark}).
Note the different behavior of the limit dispersion curves depending on $H$;   see   Figures~\ref{fig_variosH}, \ref{fig_new} and \ref{fig2} to realize the difference for the following cases:  $$ \displaystyle
a) \quad H\in\Big({ 0,  \frac{1}{\sqrt{8}}} \Big),
\qquad b)\quad H\in\Big[{ \frac{1}{\sqrt{8}}},\frac{1}{2}\Big),\qquad
c)\quad H\in\Big[\frac{1}{2},1\Big), \qquad
d)\quad H\in[1,+\infty).
$$
As a matter of fact, limit eigenvalues of multiplicity $3$ appear first for    larger values of $H$, the first one being $\Lambda=4\pi^2$ when  $H=1/2$, see Figures~\ref{fig_variosH}  and \ref{fig_new}.  This is the main reason why our final result  (cf. Theorem \ref{teoDc})  provides uniform bounds for discrepancies, in different intervals, depending on $H$,  and we focus our attention on the first eigenvalues.

\medskip

Let us describe the structure of the paper.

In Section \ref{sec2},  we  provide some preliminary results valid for any geometry of the holes, of smooth boundaries. First, we obtain  upper bounds for the eigenvalues of the perturbed problem \eqref{(11)}-\eqref{(14)}, namely,
$$
\Lambda^\varepsilon_m(\eta)\leq \Lambda_m^0(\eta)+ c_m \ee\,\,\quad \mbox{\rm for}\quad
\varepsilon \leq \varepsilon_m, \quad \eta \in [-\pi,\pi],
$$
that show how  they are controlled by those of the homogenized problem \eqref{(29)}-\eqref{(25)}.  Then, in Section \ref{sec2.1}, we introduce  a result of convergence of the eigenvalues and corresponding eigenfunctions towards those of the homogenized problem with conservation of the multiplicity,  which also  allows a perturbation of the Floquet-parameter (see Theorem \ref{corollary_convergence} and Corollary \ref{Lemma_convergence}),
\begin{equation}\label{contra}
 \Lambda^{\ee_l}_{m}(\eta_l)\to \Lambda^{0}_m (\widehat \eta),\quad \mbox{ as } \quad  (\ee_l,\eta_l)\lil (0, \widehat  \eta).
\end{equation}
This shows a  strong stability of both problems on the parameter $\eta$ which arises in the quasi-periodicity conditions, cf.  \eqref{(12-13)}  and \eqref{(25)}. The result has recently been introduced in the literature of model problems for waveguides, cf. \cite{nuestroNeumann} and \cite{GoNaOrPe_IMSE}, and  proves to be essential for obtaining uniform bounds (in the perturbation and  Floquet parameters) for convergence rates between the eigenvalues of the perturbation and limit problems. In Section \ref{sec3.1}, we introduce a boundary layer problem posed in an unbounded strip  (cf. \eqref{(33)}-\eqref{(40)}); it has been obtained by  means of asymptotic expansions in \cite{nuestroNeumann}.

In Sections \ref{sec3} and \ref{sec3.3},   we obtain  complementary results on the asymptotics of the eigenvalues which provide   bounds for the distance between the  dispersion curves of the limit problem and those of the homogenization problem. To do it, we use
a well-known result on {\it   almost eigenvalues and eigenfunctions} from the spectral perturbation theory,  cf. Lemma \ref{Lemma_Visik}, which implies the construction of  families of quasimodes from the solutions of the homogenized problem  and the boundary layer  problem, cf. \eqref{quasim}-\eqref{casiortho2as}. Some restrictions are performed both on the geometry of the hole and on the choice of the limiting eigenvalues (cf. \eqref{(symm)}  and \eqref{kpj}). Avoiding these restrictions implies that further terms of asymptotic expansions as well as further boundary layer functions are necessary, and the computations become excessively high.
Summarizing,    we can set:
$$  |\Lambda_{p(\ee,\eta,m)}^\varepsilon(\eta)-\Lambda_m^0(\eta)|<C_0 \varepsilon,\qquad \mbox{for }   \varepsilon\leq\varepsilon_0, \,\,  \eta \in  [-\pi,\pi],\, \mbox{ and some }\,  p(\ee,\eta,m)\geq m , $$
\noindent when $m=1$ and $H>0$, when $m=2$ and  $H\in(0,1/2)$  and when $m=3$ and $H\in (0,1/\sqrt{8})$.
The same results,  for $\eta$ ranging in smaller intervals (dependent on $m$  and $H$) are obtained in Section \ref{sec3.3} for higher values of $m$ (see Theorems \ref{any_j} and \ref{th_epsetabis}).

At this stage  we cannot assert that $p(\ee,\eta,m)= m$, and this is the aim of Section \ref{sec4}, where we combine the   results in Sections~\ref{sec2} and \ref{sec3.3} with different arguments of contradiction involving convergence \eqref{contra}. We obtain
$$  |\Lambda_{m}^\varepsilon(\eta)-\Lambda_m^0(\eta)|<C_0 \varepsilon,\qquad \mbox{with }  \varepsilon\leq \varepsilon_0 \,  \mbox{ and } \,\eta \in  [-\pi,\pi],  $$
for the same values of $m$ and $H$.
Above and throughout the paper,   $\ee_m$, $c_m$ and $C_m$, with $m\in  \mathbb{N}_0$, denote certain positive constants independent of both variables $ \ee$ and $\eta$.}

Due to the complexity of the  trusses-nodes structure of the dispersion curves  (cf. Figure \ref{fig_variosH}),  in the theorems, we have imposed restrictions on the index $m$  depending on $H$, which can affect the intervals where the Floquet-parameter ranges, but the technique can be extended to further values of $m$.

Finally, it should be emphasized that the method here developed can be applied to many problems arising in waveguide theory, cf. \cite{na453,na522,na537,BoPa,BaNaRu} among others.

\setcounter{equation}{0} \setcounter{section}{1}
\renewcommand{\theequation}{\arabic{section}.\arabic{equation}}
\section{Preliminary results: convergence and boundary layer.}\label{sec2}}
In Section \ref{sec2.1}, we state a result of convergence for the spectrum   of the $\ee$-dependent problem towards that of the homogenized one which allows   the perturbation of the Floquet parameter, cf. Theorem  \ref{corollary_convergence}. This result   extends the usual convergence   of the spectrum for each fixed   $\eta\in [-\pi,\pi]$, cf. Corollary \ref{Lemma_convergence}, and becomes necessary to derive our result in Section \ref{sec4}. In Section \ref{sec3.1},  we introduce a boundary layer problem and its solution  which proves  useful, in Section \ref{sec3}, for the construction of new approximations to groups of eigenfunctions of \eqref{(11)}-\eqref{(14)} corresponding to eigenvalues in small intervals.

First,  below, we obtain some  estimates for the eigenvalues of the perturbation problem which relate   homogenization and homogenized problems, as a consequence of the comparison of both spectra.

\begin{proposition}\label{minmaxu}
{\it For fixed $m=1,2,\cdots$,  let $\Lambda^\varepsilon_m(\eta)$  and $\Lambda^0_m(\eta)$   be the eigenvalues in the sequence  \eqref{(17)} and \eqref{(S2)} respectively. Then,    there exist positive constants $c_m$ and $\ee_m$, independent of $\ee$ and $\eta$, such that
\begin{equation}\begin{array}{c}\displaystyle
\Lambda^\varepsilon_m(\eta)\leq \Lambda_m^0(\eta)+ c_m \ee, \, \quad\mbox{\rm for}
\,\,\varepsilon\in(0,\varepsilon_m],\quad  \eta \in [-\pi,\pi].
\end{array}\label{(KK1)}\end{equation}
}
\end{proposition}
\begin{proof}
Let us apply the minimax principle. For each $m=1,2,\cdots$ and $\eta\in [-\pi,\pi]$, we write 
\begin{equation}\label{minmax}
\Lambda^\varepsilon_m(\eta)=
\min_{E_m^\varepsilon\subset H^{1,\eta}_{per}(\varpi^\varepsilon)} \max_{V\in E_m^\varepsilon, V\neq 0}
\frac{(\nabla_x V, \nabla_x V)_{\varpi^\varepsilon}}{(V,V)_{\varpi^\varepsilon}},
\end{equation}
where the minimum is computed over the set of subspaces $E_m^\varepsilon$ of $H^{1,\eta}_{per}(\varpi^\varepsilon)$
with dimension $m$.

To prove \eqref{(KK1)}, we take a particular subspace  $E_m^{\varepsilon,*}$ that we construct as the linear space
$$E_m^{\varepsilon,*}=\langle { U}^0_1(\cdot;\eta)\big\vert_{\varpi^\ee},  \cdots, { U}^0_m(\cdot;\eta)\big\vert_{\varpi^\ee}\rangle,$$
where   $ \{{ U}^0_k(\cdot;\eta)\}_{k=1}^m$ are the eigenfunctions of \eqref{(16limit)}, corresponding to the eigenvalues  $\{{ \Lambda }^0_k\}_{k=1}^m$ which are orthonormal in $L^2(\varpi^0)$.

For $V\in  E_m^{\varepsilon,*}$ without any restriction we can assume that $\Vert V;\, L^2( {\varpi^0}) \Vert^2=1$ (cf.   \cite{BaNaRu} and \cite{DESOME} for the idea in different problems). Thus, we have $V= \alpha_1{ U}^0_1(\cdot;\eta)+\cdots +\alpha_m U_m^0(\cdot;\eta)$ for some constants $\alpha_i$ such that  $\Vert V;\, L^2( {\varpi^0}) \Vert^2= \alpha_1^2 +\cdots +\alpha_m^2=1$ . Therefore, we write
$$(\nabla_x V, \nabla_x V)_{\varpi^\varepsilon}\leq (\nabla_x V, \nabla_x V)_{\varpi^0}=\alpha_1^2\Lambda_1^0(\eta)  +\cdots+\alpha_m^2\Lambda_m^0(\eta)\leq \Lambda_m^0(\eta).$$
Also, we can write
$$(  V,  V)_{\varpi^\varepsilon} =(  V,  V)_{\varpi^0}-(  V,  V)_{\varpi^0\setminus\varpi^\varepsilon} =1-(  V,  V)_{\varpi^0\setminus\varpi^\varepsilon}\geq 1-\widetilde c_m\ee, $$
for a certain constant $\widetilde c_m$ independent of $\eta$ and sufficiently small $\ee$. Indeed, to show the last inequality, we use the estimate
 $\, \|V;\,{L^2(\varpi^0 \cap \{|x_1|<\varepsilon\})\|}^2 \leq C\varepsilon \|V;\,{H^1(\varpi^0 )}\|^2\,$ $\forall V\in H^1( \varpi^0)$,
 see, e.g.,  Lemma~2.4 in \cite{Marchenkobook}, and consequently,
 $$(  V,  V)_{\varpi^0\setminus\varpi^\varepsilon} \leq c\ee\left((  \nabla V,  \nabla V)_{\varpi^0} + (    V,   V)_{\varpi^0}\right)\leq c\ee (\Lambda_m^0(\eta)+1)\leq \widetilde c_m \ee.$$

Hence, \eqref{minmax} leads us to
$$
\Lambda^\varepsilon_m(\eta)\leq
 \max_{V\in E_m^{\varepsilon,*}, V\neq 0}
\frac{(\nabla_x V, \nabla_x V)_{\varpi^\varepsilon}}{(V,V)_{\varpi^\varepsilon}}\leq \frac{\Lambda_m^0(\eta)}{1-\widetilde c_m\ee} \leq \Lambda_m^0(\eta) + c_m\ee,
$$
for some constant $c_m$ and $\ee$ small enough. This shows   \eqref{(KK1)} and    the proposition is proved.
\end{proof}


\subsection{ Convergence for eigenvalues.}\label{sec2.1}

\medskip

Based on extensions of eigenfunctions in perforated domains (see, for example, Section I.4.2 in \cite{Oleinikbook}),  the continuity on the Floquet-parameter of the  mappings  \eqref{(18)} and \eqref{(18m)}, and some contradiction  arguments,    the following result has been proved in \cite{nuestroNeumann}.

{
\begin{theorem}\label{corollary_convergence}
{\it For   each subsequence  $\{(\varepsilon_r, \eta_r)\}_{r=1}^\infty$ such that $\ee_r\to 0$ and $\eta_r\to \widehat  \eta \in [-\pi,\pi]$,  as $r\to \infty,$ the eigenvalues $\Lambda_m^{\varepsilon_r}(\eta_r)$ of problem \eqref{(11)}--\eqref{(14)} in the sequence \eqref{(17)},  when $(\ee,\eta)\equiv (\ee_r,\eta_r)$,  converge towards the eigenvalues of problem \eqref{(29)}--\eqref{(25)} for $\eta\equiv \widehat  \eta$, as $ r\to \infty $, and there is conservation of the multiplicity. Namely, for each $m=1,2,\cdots$, the convergence
 \begin{equation}
\label{cDr}
 \Lambda^{\ee_r}_{m}(\eta_r)\to \Lambda^{0}_m (\widehat \eta),\quad \mbox{ as }   r\to \infty,
\end{equation}
holds, where $\Lambda^0_m(\widehat  \eta)$ is the m-th eigenvalue in the sequence \eqref{(S2)}
of eigenvalues of \eqref{(29)}--\eqref{(25)} for $\eta\equiv \widehat  \eta$.

 In addition,   we can extract a subsequence, still denoted by $\varepsilon_r$, such that some suitable extension of the eigenfunctions $\{U_m^{\,\varepsilon_r}(\cdot;\eta_r)\}_{m=1}^\infty$ normalized in $L^2(\varpi^{\varepsilon_r})$,  $\{{\widehat U}_m^{\,\varepsilon_r}( \cdot; \widehat  \eta)\}_{m=1}^\infty$, converge in $L^2(\varpi^0)$, as $r\to \infty,$ towards the eigenfunctions   $\{U_m^{0}(\cdot; \widehat  \eta )\}_{m=1}^\infty$  of \eqref{(29)}--\eqref{(25)}   which form an orthonormal basis of $L^2(\varpi^0)$.}
\end{theorem}

Theorem \ref{corollary_convergence} shows a certain stability of the limit of the spectrum of \eqref{(11)}--\eqref{(14)}  under any perturbation of the Floquet-parameter $\eta$.

\begin{corollary}\label{Lemma_convergence}
{\it For any $\eta\in[-\pi,\pi]$, the eigenvalues $\Lambda_m^\varepsilon(\eta)$ of problem  \eqref{(16)} in the sequence \eqref{(17)} converge  towards the eigenvalues of problem \eqref{(29)}--\eqref{(25)} in the sequence \eqref{(S2)}, namely,
\begin{equation}\begin{array}{c}\displaystyle
\Lambda^\varepsilon_m(\eta)\,\,\to \Lambda^0_m(\eta)\,\,\quad \mbox{\rm
as }\,\,\varepsilon\to  0,
\end{array}\label{(27)}\end{equation}
 and there is conservation of the multiplicity. In addition, for each sequence, we can extract a subsequence, still denoted by $\varepsilon$, such that some suitable extension of the eigenfunctions {$\{U_m^{\,\varepsilon}(\cdot;\eta)\}_{m=1}^\infty$ normalized in $L^2(\varpi^\varepsilon)$,   $\{{\widehat U}_m^{\,\varepsilon}\}_{m=1}^\infty$ , converge in $L^2(\varpi^0)$, as $\ee\to 0$,} towards the eigenfunctions  $\{U_m^{0}(\cdot;   \eta )\}_{m=1}^\infty$ of \eqref{(29)}--\eqref{(25)} which form an orthonormal basis of $L^2(\varpi^0)$.}
\end{corollary}

This corollary    states   the classical spectral convergence for problem \eqref{(16)}; its  proof is an immediate  consequence of
Theorem~\ref{corollary_convergence}.

\medskip
\subsection{The boundary layer  problem. }\label{sec3.1}
As usual in two-scale boundary homogenization, near the perforation string $\omega^\varepsilon(0),\dots,\omega^\varepsilon(N-1)\subset \varpi^0 , $
 there appears a boundary layer which is described in the stretched coordinates
 \begin{equation}\begin{array}{c}
\xi=(\xi_1,\xi_2)=\varepsilon^{-1}{\color{black}(x_1,x_2-\varepsilon k H)}.
\end{array}\label{(33)}\end{equation}
Using these auxiliary coordinates,  we  introduce a boundary layer problem and its
solution in the unbounded perforated strip
\begin{equation}\label{def_Xi}
\Xi:=\{x\in   {\mathbb R} \times (0,H)\} \setminus\overline{\omega}\,.
\end{equation}
Obviously, for each $k=0,\cdots, N-1$, the change of variable \eqref{(33)} transforms $\omega^\varepsilon(k)$ into $\omega$, cf.  \eqref{(2)}.
 The proof of the  results  of this section    can be found in \cite{nuestroNeumann}; cf. also \cite{SOME} and \cite{SOME_IMSE} for the technique and further references.

Let us consider the function $W_0^1\ $  to be harmonic in $\Xi$,
 \begin{equation}\begin{array}{c}
-\Delta_\xi W_0^1(\xi)=0,\,\,\xi\in {\color{black}\Xi,}
\end{array}\label{(34)}\end{equation}
with the periodicity conditions
 \begin{equation}\begin{array}{c}\displaystyle
W_0^1(\xi_1,H)=W(\xi_1,0),\,\,\frac{\partial W_0^1}{\partial \xi_2}(\xi_1,H)=
\frac{\partial W_0^1}{\partial \xi_2}(\xi_1,0),\,\,\xi_1\in{\mathbb R},
\end{array}\label{(35)}\end{equation}
and the nonhomogeneous Neumann condition on the boundary of the hole
$\overline{\omega}$
  \begin{equation}\begin{array}{c}
 \partial_{\nu(\xi)} W^1_0(\xi)=-\partial_{\nu(\xi)}\xi_1=-\nu_1(\xi),\,\,\xi\in\partial\omega.
 \end{array}\label{(40)}\end{equation}
 Here,
$\nu(\xi)=(\nu_1(\xi),\nu_2(\xi))$ denote the outward (with respect to
$\Xi$) normal vector on $\partial\omega$.

Also,
  for convenience, we introduce here the  cut-off
functions, $\chi_\pm\in C^\infty({\mathbb R})$,
 \begin{equation}\begin{array}{c}\displaystyle
\chi_\pm(t)=1\,\,\mbox{\rm for}\,\,\pm t>2R\,\,\mbox{\rm and}\,\, \chi_\pm(t)=0\,\,\mbox{\rm for}\,\,\pm t<R,
\end{array}\label{(chi)}\end{equation}
with a fixed   $R>0$ satisfying
 \begin{equation}\begin{array}{c}\displaystyle
 \overline{\omega}\subset\Xi(R):=\{\xi\in\Xi:\,|\xi_1|<R\}.
\end{array}\label{(ksi)}\end{equation}

The variational formulation of \eqref{(34)}-\eqref{(40)}  reads: find $W_0^1\in {\cal H}^1_{per}(\Xi)$ satisfying
 \begin{equation}\begin{array}{c}\displaystyle
 (\nabla_\xi W_0^1,\nabla_\xi V)_\Xi= (-\nu_1(\xi),V)_{\partial \color{black} \omega}\quad
 \forall V\in {\cal H}^1_{per}(\Xi),
\end{array}\label{(44po)}\end{equation}
where the space ${\cal H}^1_{per}(\Xi)$  denotes the completion of the
linear space $C^\infty_{per}(\overline{\Xi})$ (of the infinitely
differentiable $H$-periodic in $\xi_2$ functions with compact
supports) in the norm
  \begin{equation}\begin{array}{c}\displaystyle
 \|W;{\cal H}^1_{per}(\Xi)\|=\Big(\|\nabla_\xi W;L^2(\Xi)\|^2+\|W;L^2(\Xi(2R))\|^2\Big)^{1/2}.
\end{array}\nonumber
\end{equation}

Since the compatibility condition, $ (-\nu_1(\xi),1)_{\partial \omega  }=0,$  is satisfied,   problem   \eqref{(44po)}
has a unique solution   $W^1_0\in{\cal H}^1_{per}(\Xi)$ which is {\color{black} uniquely} defined
up to an additive constant. Moreover, since the boundary $\partial\omega$ is
smooth, this solution is infinitely differentiable in
$\overline{\Xi}$ and the Fourier method (cf. \cite{LeguillonSP}, \cite{MaNaPl} and \cite{nuestroNeumann}), in particular, gives the
decomposition
   \begin{equation}\begin{array}{c}
 W^1_0(\xi)=\sum\limits_\pm\chi_\pm(\xi_1)C_\pm+\widetilde{W}^1_0(\xi)
   \end{array}\label{(W100)}\end{equation}
with the exponentially decaying remainder $\color{black} \widetilde{W}^1_0$, and some constants
$ C_\pm$ {\color{black} which can also depend on $R$, {\color{black}cf.~\eqref{(chi)} and \eqref{(ksi)}}.}

\begin{figure}[t]
\begin{center}
\scalebox{0.45}{\includegraphics{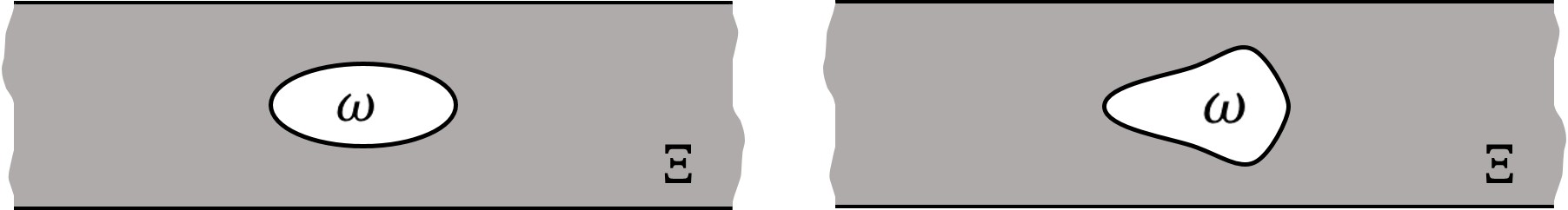}}
\end{center}
\caption{The strip $\Xi$ with two different possible geometries for the hole $\omega$.}
\label{fig-mirror}
\end{figure}

Note that the above results hold for any smooth hole. In addition,
 under the assumption of  {\it mirror  symmetry}:
 \begin{equation}\begin{array}{c}\displaystyle
\omega=\{\xi=(\xi_1,\xi_2)\in{\mathbb
R}^2:\,(\xi_1,H-\xi_2)\in\omega\},
\end{array}\label{(symm)}\end{equation}
 see, for instance,  Figure~\ref{fig-mirror}, the
function $W_0^1$ is even  in the
 $\xi_2-H/2$  variable and satisfies
\begin{equation} \displaystyle
\frac{\partial W_0^1}{\partial \xi_2}(\xi_1,0)=\frac{\partial
W_0^1}{\partial \xi_2}(\xi_1,H)=0, {\color{black}\,\,\, \xi_1\in\mathbb{R},}
 \label{(J8)}\end{equation}
 see \cite{nuestroNeumann} for a proof.

\setcounter{equation}{0} \setcounter{section}{2}
\renewcommand{\theequation}{\arabic{section}.\arabic{equation}}

\section{Asymptotic formulas for convergence rates}\label{sec3}

In this and  next section,  we  obtain some   important  complementary results of \eqref{(27)}  on  the asymptotics of   the perturbed  dispersion
curves, at a low  frequency range. This implies constructing  the so-called   {\it almost eigenvalues and eigenfunctions}  or {\it quasimodes} of a compact operator associated with \eqref{(11)}-\eqref{(14)}.   In order to do this, we set an abstract framework for the perturbation problem and construct the approximations from the eigenfunction of the homogenized problem and a boundary layer function (cf. Section \ref{sec3.1}).    We obtain asymptotic formulas for eigenvalues and estimates when the limit eigenvalues have corresponding eigenfunctions which do not depend on $x_2$, cf. \eqref{(31)},  \eqref{kpj} and \eqref{(un3)}. These formulas   do not allow  us to predict the precise index of the eigenvalue in the sequence \eqref{(17)}, which  will be determined in Section \ref{sec4} using the results of Section \ref{sec2}.
Another restriction that simplifies  computations in this section comes from the geometry of the perforation:   we assume the mirror symmetry (cf. Figure \ref{fig-mirror}, \eqref{(symm)} and \eqref{(J8)}).

First, we   reformulate the spectral problem \eqref{(11)}--\eqref{(14)} in terms of  operators on Hilbert spaces.  Let ${\cal
H}^\varepsilon(\eta)$  denote the  space  $H^{1,\eta}_{per}(\varpi^\varepsilon)$ equipped with the
usual scalar product in $H^1(\varpi^\ee)$. We denote by $\langle\cdot, \cdot\rangle_{\eta\ee}$ this scalar product and by  $\|U^\varepsilon;{\cal
H}^\varepsilon(\eta)\|$    the associated norm.

In the space ${\cal H}^\varepsilon(\eta)$, we define the compact, positive and  symmetric
operator ${\cal B}^\varepsilon(\eta)$, as
\begin{equation}\begin{array}{c}\displaystyle
\langle {\cal B}^\varepsilon(\eta)U^\varepsilon,V^\varepsilon\rangle_{\eta\ee}=
( U^\varepsilon,V^\varepsilon)_{\varpi^\varepsilon}\quad\forall U^\varepsilon,V^\varepsilon
\in H^{1,\eta}_{per}(\varpi^\varepsilon).
\end{array}\label{(J2)}\end{equation}
It is self-evident that the variational formulation \eqref{(16)} of \eqref{(11)}--\eqref{(14)}
can be re-written as follows:
\begin{equation}\begin{array}{c}\displaystyle
{\cal B}^\varepsilon(\eta)U^\varepsilon(\eta)=M^\varepsilon(\eta)
U^\varepsilon(\eta)\,\, \quad \mbox{\rm in} \quad \,\,{\cal H}^\varepsilon(\eta),
\end{array}\label{(J3)}\end{equation}
with the new spectral parameter
\begin{equation}\begin{array}{c}\displaystyle
M^\varepsilon(\eta)=(1+\Lambda^\varepsilon(\eta))^{-1}.
\end{array}\label{(J4)}\end{equation}

The following result, which we state for the specific operator ${\cal B}^\varepsilon(\eta)$ in \eqref{(J2)}, is  based on the lemma on almost eigenvalues and eigenfuntions from the spectral perturbation theory, cf. \cite{ViLu},  \cite[\S 3.1]{Oleinikbook}, \cite[\S 7.1]{Nazarovbook} and \cite[\S 5.32]{Lazutkin}.

\medskip

\begin{lemma}\label{Lemma_Visik}
{\it Let $M^\varepsilon_{as}(\eta)\in{\mathbb R}$
and $U^\varepsilon_{as}(\eta) \in{\cal
H}^\varepsilon(\eta)\setminus\{0\}$ verify the relationship
\begin{equation}\begin{array}{c}\displaystyle
\|{\cal
B}^\varepsilon(\eta)U^\varepsilon_{as}(\eta)-M^\varepsilon_{as}(\eta)
U^\varepsilon_{as}(\eta);{\cal H}^\varepsilon(\eta)\|=\delta^\ee\|
U^\varepsilon_{as}(\eta);{\cal H}^\varepsilon(\eta)\|.
\end{array}\label{(J5)}\end{equation}
Then, there exists an eigenvalue $M^\varepsilon(\eta)$ of the
operator ${\cal B}^\varepsilon(\eta)$ such that
\begin{equation}\label{J51} \displaystyle
|M^\varepsilon(\eta)-M^\varepsilon_{as}(\eta)|\leq\delta^\ee.
\end{equation}

Moreover, for any $\delta_1^\ee>\delta^\ee$, there exist coefficients $\alpha_J^\ee,\cdots, \alpha_{J+K-1}^\ee$, satisfying
\begin{equation}\begin{array}{c}\displaystyle
\left\|
\|U^\varepsilon_{as}(\eta);{\cal H}^\varepsilon(\eta)\|^{-1}U^\varepsilon_{as}(\cdot;\eta)-\sum_{j=J}^{J+K-1} \alpha_j^\ee U_j^\varepsilon (\eta)
 ;{\cal H}^\varepsilon(\eta)\right\|\leq 2\frac{\delta^\ee}{\delta_1^\ee}, \quad \sum_{j=J}^{J+K-1}(\alpha_j^\ee)^2 =1,
\end{array}\label{(J5e)}\end{equation}
where  $M_J^\varepsilon(\eta) ,\cdots,  M_{J+K-1}^\varepsilon(\eta)$
are all the eigenvalues of the operator ${\cal B}^\ee (\eta)$    in the interval $[M_{as}^\varepsilon(\eta)-\delta_1^\ee, M_{as}^\varepsilon(\eta)+\delta_1^\ee]$, and  $U_J^\varepsilon(\cdot; \eta),  \cdots , U_{J+K-1}^\varepsilon(\cdot;\eta)$ are the corresponding eigenvectors normalized by
$$\langle U_i^\varepsilon(\cdot;\eta),U_j^\varepsilon(\cdot; \eta)\rangle_{\eta\ee}=\delta_{ij}, \quad i,j=J,  \cdots, J+{K-1},$$
being $\delta_{ij}$  the Kronecker delta.}
\end{lemma}

The pair   $(\, M^\varepsilon_{as}(\eta),\,  \|
U^\varepsilon_{as}(\eta);{\cal H}^\varepsilon(\eta)\|^{-1}U^\varepsilon_{as}(\eta)\, )$  in  Lemma \ref{Lemma_Visik}  is    the  so-called quasimode   of operator ${\cal B}^\ee(\eta)$ with remainder $\delta^\ee$. It approaches    eigenvalues and eigenfunctions of the operator  ${\cal B}^\ee(\eta)$ as stated in \eqref{J51}-\eqref{(J5e)}. Also, if there is no confusion, the function $ \|
U^\varepsilon_{as}(\eta);{\cal H}^\varepsilon(\eta)\|^{-1}U^\varepsilon_{as}(\eta)$ is referred to as the quasimode.

In what follows,  throughout the section,  we construct  $M^\varepsilon_{as}(\eta)$ and
$U^\varepsilon_{as}(\eta)$ and obtain a bound for $\delta^\ee$ in \eqref{(J5)}, cf. \eqref{deltacota}. The relation of $M^\varepsilon_{as}(\eta) $ with the true eigenvalues in the sequence \eqref{(17)} will be given   in Section \ref{sec3.3}.

\bigskip

For $ j=0,1,2,\cdots$, let
\begin{equation}\label{kpj}\Lambda^0_{\pm,j}(\eta)=({\color{black}\eta}\pm 2\pi j)^2\end{equation} be   eigenvalues  in
\eqref{(31)} corresponding to a fixed Floquet parameter
$\eta\in[-\pi ,\pi]$.  From the explicit computations, the corresponding eigenfunctions
\begin{equation}\begin{array}{c}\displaystyle
U^0_{\pm,j}(x_1;\eta)=e^{{\rm i}(\eta\pm  {2}\pi j)x_1},
\end{array}\label{(un3)}\end{equation}
do not depend on the $x_2$-variable. Obviously, for the different signs, values in \eqref{kpj} and functions in \eqref{(un3)}, coincide only for $j=0$. Namely, for $j=0$, we have the first eigenvalue in the increasing sequence \eqref{(S2)} and the corresponding eigenfunction:
\begin{equation}\label{kp0} \Lambda_1^0(\eta)= \Lambda^0_{\pm,0}(\eta)= \eta  ^2,  \quad U_1^0(\cdot; \eta)   =U^0_{\pm,0}(x_1;\eta)  =e^{{\rm i} \eta x_1}.\end{equation}

 According to the notation in \eqref{(J4)} we take
\begin{equation}\begin{array}{c}\displaystyle
M_{\pm,j}^0(\eta)=(1+\Lambda^0_{\pm,j}(\eta))^{-1}
\end{array}\nonumber\end{equation}
as an approximate eigenvalue \textcolor{black}{($M_{as}^\ee(\eta)$ for $\pm$ respectively),} and
\begin{equation}\label{(un2)}
\begin{split}
U^\varepsilon_{\pm,j}(x;\eta)=&X^\varepsilon(x_1)
U^0_{\pm,j}(x_1;\eta)+ (1-X^\varepsilon(x_1))\bigg(U^0_{\pm,j}(0;\eta)+x_1
\frac{\partial U^0_{\pm,j}}{\partial x_1}(0;\eta)\bigg) \\
&+\varepsilon\chi_0(x_1)\frac{\partial U^0_{\pm,j}}{\partial x_1}(0;\eta)
W^1_0\Big(\frac{x}{\varepsilon}\Big),
\end{split}
\end{equation}
as an approximate eigenfunction $U_{as}^\ee(\eta)$ which we construct from $U^0_{\pm,j}$ in \eqref{(un3)}  and  $W^1_0(\xi)$,
the bounded harmonics in $\Xi$, see \eqref{def_Xi}, \eqref{(34)}-\eqref{(40)} and \eqref{(W100)}. Here,
\begin{equation}
\begin{split}
& X^\varepsilon(x_1)=1-\chi_+(x_1/\varepsilon)-\chi_-(x_1/\varepsilon), \quad  \mbox{ and, } \quad \\
&\chi_0\in C^\infty({\mathbb R}),\,\, \chi_0(x_1)=1\,\,\mbox{\rm
for}\,\,|x_1|\leq 1/6,\,\, \chi_0(x_1)=0\,\,\mbox{\rm
for}\,\,|x_1|\geq 1/3,
\end{split}
\label{(ra3)}\end{equation}
where the even smooth cut-off functions $\chi_{\pm}$ are defined by \eqref{(chi)}. It can be easily verified that $U^\varepsilon_{\pm,j}(x;\eta)\in {\cal
H}^\varepsilon(\eta)$.

Note  that,  depending on $H$,  $j\in  {\mathbb N}_0$ and $\eta\in [-\pi, \pi]$, the  eigenvalues    $M_{\pm,j}^0(\eta)$ can be simple, or   have a multiplicity greater than or equal to $2$.  Obviously, once $H$ and $j$ are fixed, also the eigenvalue number in the sequence \eqref{(S2)} may change depending on $\eta$. Below, we fix $j$ and the sign plus or minus, and for brevity,  we omit to write the $\eta$ dependence of function $U^\ee_{\pm,j}$.

\medskip

 In order to apply Lemma~\ref{Lemma_Visik},   we  multiply \eqref{(J5)}  by $ \|
U^\varepsilon_{as}(\eta);{\cal H}^\varepsilon(\eta)\|^{-1}$ , write $ \delta^\ee\equiv\delta^\varepsilon_{\pm,j}(\eta)$ and    obtain the remainder
 \begin{eqnarray}\label{(un4)}
\delta^\varepsilon_{\pm,j}(\eta)&:=&\|U_{\pm,j}^\varepsilon;{\cal
H}^\varepsilon(\eta)\|^{-1}\|{\cal
B}^\varepsilon(\eta)U_{\pm,j}^\varepsilon-M^0_{\pm,j}(\eta)
U_{\pm,j}^\varepsilon;{\cal H}^\varepsilon(\eta)\|,
\end{eqnarray}
for which we obtain the  uniform estimate: \begin{equation}\label{deltacota} \delta^\varepsilon_{\pm,j}(\eta)\leq c_j\varepsilon, \quad \mbox{ for } \ee\leq \ee_j,   \end{equation} where  $c_j $  and $\ee_j$ are two    positive constants  independent of $\eta$ and $\ee$.

To prove \eqref{deltacota}, first, let us show the  almost orthogonality  property for the family of functions constructed in \eqref{(un2)}:
\begin{eqnarray}
&|\langle U_{+,j}^\varepsilon,U_{-,j}^\varepsilon\rangle_{\eta\ee}|\leq  C_j\varepsilon^{1/2}\qquad \mbox{for }\varepsilon<\varepsilon_j, \, \eta\in [-\pi,\pi], \quad j \in \mathbb{N} ,
\label{casiortho} \vspace{0.3cm} \\
&|\langle U_{ \pm,j}^\varepsilon,U_{ \mp ,k}^\varepsilon\rangle_{\eta\ee}|\leq  C_{j,k}\varepsilon^{1/2}\qquad \mbox{for }\varepsilon<\varepsilon_{j,k}, \, \eta\in [-\pi,\pi], \quad j,k\in \mathbb{N}_0,\, j\not=k, \label{casiortho2}
\end{eqnarray}
where $C_j$, $C_{j,k}$, $\ee_j$ and $\ee_{j,k}$ are some
  positive  constants independent of $\varepsilon$ and $\eta $; recall that $\langle \cdot,\cdot\rangle_{\eta\ee}$ denotes the scalar product in  $\mathcal{H}^\varepsilon(\eta)$.

Let us prove \eqref{casiortho}. Owing to the orthogonality of the functions $U_{\pm,j}^0$ in $L^2(\varpi^0)$ and $H^1(\varpi^0)$, we write
\begin{equation}
\begin{split}
|\langle U_{+,j}^\varepsilon,U_{-,j}^\varepsilon\rangle_{\eta\ee}|= &(\nabla_x(U_{+,j}^\varepsilon-U_{+,j}^0), \nabla_x U_{-,j}^\varepsilon)_{\varpi^\varepsilon}+
(\nabla_x U_{+,j}^0, \nabla_x (U_{-,j}^\varepsilon-U_{-,j}^0))_{\varpi^\varepsilon}   \\&- (\nabla_x U_{+,j}^0, \nabla_x U_{-,j}^0)_{\varpi^0\setminus\varpi^\varepsilon}
+(U_{+,j}^\varepsilon-U_{+,j}^0, U_{-,j}^\varepsilon)_{\varpi^\varepsilon} \\&+
(U_{+,j}^0, U_{-,j}^\varepsilon-U_{-,j}^0)_{\varpi^\varepsilon} - (U_{+,j}^0,U_{-,j}^0)_{\varpi^0\setminus\varpi^\varepsilon}. \label{decomp}
\end{split}
\end{equation}
In addition, by the definition \eqref{(un2)} of $U_{\pm,j}^\varepsilon$, we have
\begin{equation*}
\begin{split}
\|U_{\pm,j}^\varepsilon-U_{\pm,j}^0; \mathcal{H}^\varepsilon(\eta)\|\leq &
\bigg\|(1-X^\varepsilon(x_1))\bigg(U^0_{\pm,j}(x_1;\eta)-U^0_{\pm,j}(0;\eta)-x_1
\frac{\partial U^0_{\pm,j}}{\partial x_1}(0;\eta)\bigg); \mathcal{H}^\varepsilon(\eta)\bigg\| \\
&+\bigg\|\varepsilon\frac{\partial U^0_{\pm,j}}{\partial
x_1}(0;\eta)W^1_0\big(\frac{x}{\varepsilon}\big) ; \mathcal{H}^\varepsilon(\eta)\bigg\|.
\end{split}
\end{equation*}
Thus, since   $U_{\pm,j}^0$ are smooth functions,  the support of $1-X^\varepsilon$ is contained in $\{|x_1|<2R\varepsilon\}$ and $W_0^1\in \mathcal{H}_{per}^1(\Xi)$, we have \begin{equation}\label{dif}
\|U_{\pm,j}^\varepsilon-U_{\pm,j}^0; \mathcal{H}^\varepsilon(\eta)\|\leq \widehat C_j\varepsilon^{1/2},
\end{equation}
where,  on account of  \eqref{(ra3)},  we have applied the Taylor formula for $U_{\pm,j}^0(x_1,\eta)$ while $\vert x_1\vert\leq O(\ee)$, the change of variable \eqref{(33)}, and the periodicity of $W_0^1(\xi)$ in the $\xi_2$-direction. Also $ \widehat C_j$ is a constant independent of $\ee$ and $\eta$.
Now, using  \eqref{dif}, the smoothness of $U_{\pm,j}^0(x_1,\eta)$ and  the   fact that $|\varpi^0\setminus\varpi^\varepsilon|=O(\varepsilon)$, we obtain
\begin{equation}\label{(un64)}
\|U_{\pm,j}^\varepsilon;{\cal H}^\varepsilon(\eta)\|^2 \lie \|U_{\pm,j}^0;L^2(\varpi^0)\|^2+ \|\nabla_x U_{\pm,j}^0;L^2(\varpi^0)\|^2=(1+\Lambda_{\pm,j}^0(\eta))H.
\end{equation}
Hence, gathering \eqref{decomp}, \eqref{dif}, \eqref{(un64)} and using that $|\varpi^0\setminus\varpi^\varepsilon|=O(\varepsilon)$, we get \eqref{casiortho}.

Rewriting the proof above, with minor modifications, for each $k$ and $j$, $\k\not= j$,  we obtain the four estimates in \eqref{casiortho2}.

Then, for each sign plus or minus, index $j\in \mathbb{N}_0$ and $\eta\in[-\pi,\pi]$,  let us introduce
\begin{equation}\label{quasim}U_{as,\pm,j}^\ee(\eta):= \|U_{\pm,j}^\varepsilon(\eta);{\cal H}^\varepsilon(\eta)\|^{-1} U^\ee_{\pm,j}(\eta) \end{equation}
to be the quasimode constructed from the eigenfunction corresponding with  $\Lambda_{\pm,j}^0(\eta)$, cf. \eqref{(un2)} and \eqref{(un3)}. From \eqref{casiortho},  \eqref{casiortho2} and \eqref{(un64)},  we get the    almost orthonormality  conditions:
 \begin{eqnarray}
&|\langle U_{as.+,j}^\varepsilon,U_{as,-,j}^\varepsilon\rangle_{\eta\ee}|\leq \widetilde C_j\varepsilon^{1/2}\qquad \mbox{for }\varepsilon<\widetilde \varepsilon_j, \, \,\eta\in [-\pi,\pi],\quad j\in \mathbb{N}, \label{casiorthoas} \vspace{0.3cm}\\
&|\langle U_{as, \pm,j}^\varepsilon,U_{as, \mp ,k}^\varepsilon\rangle_{\eta\ee}-\delta_{kj}|\leq  \widetilde C_{j,k} \varepsilon^{1/2}\qquad \mbox{for }\varepsilon<\widetilde \varepsilon_{j,k}, \, \,\eta\in [-\pi,\pi], \quad j,k\in \mathbb{N}_0, \label{casiortho2as}
\end{eqnarray}
where $\widetilde C_j$, $\widetilde C_{j,k}$, $\widetilde \ee_j$ and $\widetilde \ee_{j,k}$ are some
positive  constants independent of $\varepsilon$ and $\eta$.

Finally, we write the following estimate
$$  \|{\cal
B}^\varepsilon(\eta)U_{\pm,j}^\varepsilon-M^0_{\pm,j}(\eta)
U_{\pm,j}^\varepsilon;{\cal H}^\varepsilon(\eta)\|\leq \widehat c_j\ee,   \quad \mbox{ for } \ee\leq \widehat \ee_j $$ (for a certain positive constant  $\widehat c_j$ independent of $\varepsilon$ and $\eta$), whose proof involves
cumbersome computations that we avoid introducing them here: it follows rewriting the proof in \cite{nuestroNeumann} with minor modifications. This estimate together with  the convergence \eqref{(un64)} allows  us to write the uniform bound for the remainder \eqref{deltacota}.

\medskip

Now,  applying   Lemma~\ref{Lemma_Visik} gives us an eigenvalue
$M^\varepsilon_{\pm,j}(\eta)$ of the operator ${\cal
B}^\varepsilon(\eta)$ such that
\begin{equation}\begin{array}{c}\displaystyle
|M^\varepsilon_{\pm,j}(\eta)-M^0_{\pm,j}(\eta)|\leq c_j\varepsilon, \quad   \mbox{ for }\ee\leq \ee_j, \, \, \eta\in [-\pi,\pi],
\end{array}\label{(un7)}\end{equation}
where the factor $c_j>0$ and the constant $\ee_j>0$ are independent of $\eta$ and $\ee$. Recalling
\eqref{(J4)}, we derive from \eqref{(un7)} that
\begin{equation} \displaystyle
|\Lambda^\varepsilon_{\pm,j}(\eta)-\Lambda^0_{\pm,j}(\eta)|\leq
c_j\varepsilon(1+
\Lambda^0_{\pm,j}(\eta))(1+\Lambda^\varepsilon_{\pm,j}(\eta)),
\label{(un8)}\end{equation}
{\color{black}  and, hence
$$\displaystyle
(1+\Lambda^\varepsilon_{\pm,j}(\eta))(1-c_j\varepsilon(1+
\Lambda^0_{\pm,j}(\eta))\leq 1+\Lambda^0_{\pm,j}(\eta).$$
Let us set
\begin{equation}\nonumber
\widetilde \varepsilon_j:=\frac{1}{2c_j(1+4\pi^2j^2)} \mbox{ when }  j=1,2,\cdots \quad \mbox { while } \widetilde \varepsilon_0:=\frac{1}{2c_0(1+ \pi^2 )} \mbox{ when }  j=0.
\label{(un9)}\end{equation}
Then, for $\varepsilon\leq \min (\ee_j,\widetilde \ee_j) \ $, we have $(1-c_j\varepsilon(1+\Lambda^0_{\pm,j}(\eta))>1/2$ and therefore
\begin{equation}\begin{array}{c}\displaystyle
|\Lambda^\varepsilon_{\pm,j}(\eta)-\Lambda^0_{\pm,j}(\eta)|\leq
{\color{black}2}c_j\varepsilon(1+ \Lambda^0_{\pm,j}(\eta))^2 \leq  {\color{black}C_j
}\varepsilon
\qquad \forall\eta\in [-\pi,\pi].
\end{array}\label{(un10)}\end{equation}

\medskip
Estimates \eqref{(un7)} and \eqref{(un10)} become essential in what follows  to derive bounds for convergence rates between the eigenvalues in sequences \eqref{(17)} and \eqref{(S2)}.

\setcounter{equation}{0} \setcounter{section}{3}
\renewcommand{\theequation}{\arabic{section}.\arabic{equation}}
\section{The almost eigenvalues and the  true eigenvalues}\label{sec3.3}
In this section,  we provide a certain relation between the eigenvalue  $\Lambda^\varepsilon_{\pm,j}(\eta)$ given by Lemma \ref{Lemma_Visik}, which satisfies  \eqref{(un10)}, and an eigenvalue in the increasing sequence \eqref{(17)}.   To state Theorems \ref{any_j} and \ref{th_epsetabis} and Corollary \ref{cor_doble},  we use  the results in Section \ref{sec3}    and arguments on families of almost orthogonal functions approaching the eigenfunctions,  cf.   \cite{LoPeCR}, \cite{PeDCDS} and \cite{SOME} for the technique in very different problems. See also     \cite[\S 5.32]{Lazutkin} in this connection.

Based on definitions \eqref{kpj}-\eqref{kp0} and estimates \eqref{(un7)} and \eqref{(un10)}, we show the following result.

\begin{theorem}\label{any_j}
{\it For each eigenvalue  $\Lambda_m^0(\eta)$  of problem  \eqref{(29)}--\eqref{(25)} in the sequence
\eqref{(S2)}, with $m=1,2,3,\cdots,$ and such that a corresponding eigenfunction does not depend of  the $x_2$-variable  for $\eta$ in an interval $[\mathfrak{a},\mathfrak{b}]\subset[-\pi,\pi]$, there is at least one  eigenvalue $\Lambda_{p(\ee,\eta,m)}^\varepsilon(\eta)$ of  problem \eqref{(11)}--\eqref{(14)}     satisfying
\begin{equation}\label{conm_j} |\Lambda_{p(\ee,\eta,m)}^\varepsilon(\eta)-\Lambda_m^0(\eta)|<C_m \varepsilon\qquad \forall    \, \, \,  \eta \in  [\mathfrak{a},\mathfrak{b}],\,  \varepsilon\leq\varepsilon_m,\end{equation}
where $C_m$ and $\ee_m$ are certain positive constants that are independent of $\ee$ and $\eta$. Moreover, for intervals  $ [\mathfrak{a},\mathfrak{b}]$ which does  not contain $\eta_0$ such that the multiplicity of $\Lambda_m^0(\eta_0)$ is greater than $1$, $p(\ee,\eta,m)\geq m$. }
\end{theorem}

\begin{proof}
Since by hypothesis, $\Lambda_m^0(\eta)$ coincides with some $\Lambda^0_{\pm,j}(\eta)$ in $[\mathfrak{a},\mathfrak{b}]$ for some sign plus or minus and some $j\equiv j(m)$, we set $\Lambda_{p(\ee,\eta,m)}^\varepsilon(\eta)=\Lambda^\ee_{\pm,j}(\eta)$, for the same sign and $j$, where it should be noted that the sign could change in subintervals of $[\mathfrak{a},\mathfrak{b}]$. Therefore, \eqref{(un10)} provides \eqref{conm_j}. The   the fact that $p(\ee,\eta,m)$ be greater than or equal to $m$ is due to the estimate \eqref{(KK1)}, and the theorem is proved. \end{proof}

It should be emphasized that while Corollary \ref{Lemma_convergence} provides approaches between all the eigenvalues of \eqref{(11)}-\eqref{(14)} and \eqref{(29)}-\eqref{(25)},  Theorem \ref{any_j} only provides estimates for discrepancies between  certain  eigenvalues of  the homogenization problem \eqref{(11)}-\eqref{(14)} and certain eigenvalues of the homogenized problem  \eqref{(29)}-\eqref{(25)} for $\eta$ in certain intervals which  can coincide with $[-\pi,\pi]$ depending on $m$ and $H$ (see Figures \ref{fig_variosH}, \ref{fig_new} and \ref{fig2}). Nevertheless, we cannot assure that all the eigenvalues \eqref{(S2)} enjoy of such an approach. In addition, we cannot assure yet that the number $p(\ee,\eta,m)$ coincides with $m$. To show this rigourously  can only be done for a few values of $m$, always depending on the value of $H$. This is due to the difficulty in ordering the eigenvalues  \eqref{kpj}  in the monotone sequence \eqref{(S2)}.  In particular, to outline the difficulty, we show that  $p(\ee,\eta,m)=m$  for the values of $m=1,2,3$ (cf. Theorem \ref{teoDc}). Next  theorem provides the preliminary  result $p(\ee,\eta,m)\geq m$ for $\eta$ in a neighborhood of $\eta_0$ such that  $\Lambda_m^0(\eta_0)$ is a multiple eigenvalue of \eqref{(29)}-\eqref{(25)}. However, on account of \eqref{casiorthoas}-\eqref{(un10)}, the process can be continued for the values of $m$ arising in the statement of  Theorem \ref{any_j}, provided that the $ \{\Lambda_{p}(\eta)\}_{p=1}^m$ has a corresponding  eigenfunction depending only on $x_1$, see, e.g.,  Figure \ref{fig_new} when $H<1/\sqrt{8}$, and Remark \ref{rmark}.

For the sake of simplicity in the proof of the next theorem, while $m=1,2,3$, it proves useful to write \eqref{(un7)} and \eqref{(un10)} as
\begin{equation}\begin{array}{c}\displaystyle
|M^\varepsilon_{\pm,j}(\eta)-M^0_{\pm,j}(\eta)|\leq
c\varepsilon
\qquad \forall\eta\in [-\pi,\pi],\, \ee\leq \ee_0
\end{array}\label{(un70)}\end{equation}
and
\begin{equation}\begin{array}{c}\displaystyle
|\Lambda^\varepsilon_{\pm,j}(\eta)-\Lambda^0_{\pm,j}(\eta)|\leq
  {\color{black}C
}\varepsilon
\qquad \forall\eta\in [-\pi,\pi],\, \ee\leq \ee_0,
\end{array}\label{(un100)}\end{equation}
for certain positive constants $\ee_0$, $c$ and $C$.
Consequently, for any $c_r\geq c$ the interval   $[M^0_{\pm,j}(\eta)-c_r\ee, M^0_{\pm,j}(\eta)+c_r\ee]$ contains at least one eigenvalue $M^\ee_{\pm,j}(\eta)$ and the value of  $\ee_0 $ can be replaced by $\ee_{0,r}$ in order that \eqref{(un100)} be satisfied for a certain constant $C_r$. Namely,
\begin{equation}\begin{array}{c}\displaystyle
|M^\varepsilon_{\pm,j}(\eta)-M^0_{\pm,j}(\eta)|\leq
c_r\varepsilon
\qquad \forall\eta\in [-\pi,\pi],\, \ee\leq \ee_{0,r}
\end{array}\label{(un70r)}\end{equation}
and
\begin{equation}\begin{array}{c}\displaystyle
|\Lambda^\varepsilon_{\pm,j}(\eta)-\Lambda^0_{\pm,j}(\eta)|\leq
  {\color{black}C_r
}\varepsilon
\qquad \forall\eta\in [-\pi,\pi],\, \ee\leq \ee_{0,r}.
\end{array}\label{(un10r)}\end{equation}
Note that the above formulas contain  the case where $j=0$, cf. \eqref{kp0}:
\begin{equation}\label{(para1)}
|(1+\Lambda^\varepsilon_{\pm,0}(\eta))^{-1}-(1+\eta^2)^{-1}|\leq
c_r\varepsilon
\qquad \forall\eta\in [-\pi,\pi],\, \ee\leq \ee_{0,r},
\end{equation}
and
\begin{equation}\label{(un10bis)}
|\Lambda^\ee_{\pm,0} (\eta)-\eta^2|\leq C_r \varepsilon \qquad \forall\eta\in [-\pi,\pi], \,\, \varepsilon \leq \varepsilon_{0,r}.
\end{equation}
}

 Similarly, for simplicity, for a fixed $m=1,2,3$ and   $\eta \in  [\mathfrak{a},\mathfrak{b}]\subset [-\pi,\pi]$,  we avoid writing  corresponding signs and index $j\equiv j(m)$ and  we denote by $U_{as,m}^\ee(\eta )  $ and $U_{as,m+1}^\ee(\eta )  $ the two  quasimodes in \eqref{quasim} constructed
from the two eigenfunction corresponding to  $\Lambda_m^0(\eta )$ and $\Lambda_{m+1}^0(\eta )$, cf. \eqref{(un2)} and \eqref{(un3)}, for the associated $j$ and sign plus or minus. On account of \eqref{casiorthoas} and \eqref{casiortho2as}   they satisfy
\begin{equation}\label{(un64as)}
\langle U_{as,m}^\varepsilon(\eta)\,,\, U_{as,m+1}^\varepsilon(\eta) \rangle_{\eta\ee}= \widetilde C  \sqrt{\ee}\to 0  \quad \mbox{ and } \quad  \|U_{as,m}^\ee(\eta); {\cal H}^\varepsilon(\eta)\|^2 =1, \quad m=1,2,3.
\end{equation}

{

\begin{theorem}\label{th_epsetabis}
{\it There exist constants $\varepsilon_0>0$ and $C_0$  independent of $\varepsilon$ and $\eta \in [\mathfrak{a},\mathfrak{b}]\subseteq [-\pi,\pi]$  such that for each  eigenvalue  $\Lambda_m^0(\eta)$  of problem  \eqref{(29)}--\eqref{(25)} in the sequence
\eqref{(S2)}, with $m=1,2,3  ,$ and a corresponding eigenfunction depending only on the $x_1$-variable, there is an eigenvalue $\Lambda_{p(\ee,\eta,m)}^\varepsilon(\eta)$ of  problem \eqref{(11)}--\eqref{(14)}     satisfying
\begin{equation}\label{conm} |\Lambda_{p(\ee,\eta,m)}^\varepsilon(\eta)-\Lambda_m^0(\eta)|<C_0 \varepsilon\qquad \mbox{for } m=1,2,3,   \, \, 0<\varepsilon<\varepsilon_0 \, \mbox{ and } \, \,\eta  \in [\mathfrak{a},\mathfrak{b}]\subset [-\pi,\pi],\end{equation}
where  $p(\ee,\eta,m)\geq m$.
Depending on $H$ and $m$,  the interval $[\mathfrak{a},\mathfrak{b}]$ can coincide with the whole $[-\pi,\pi]$ or with any interval which does not contain abscises of intersecting points of the dispersion curves with one of the corresponding eigenfunctions depending on $x_2$.

  In particular, for $m=1$ and $H>0$,  for $m=2$ and  $H\in(0,1/2)$  and for $m=3$ and $H\in (0,1/\sqrt{8})$ we get
\begin{equation}\label{conm12} |\Lambda_{p(\ee,\eta,m)}^\varepsilon(\eta)-\Lambda_m^0(\eta)|<C_0 \varepsilon,\qquad \mbox{for }  0<\varepsilon<\varepsilon_0, \,  \eta \in  [-\pi,\pi],\, \mbox{ and }\,  p(\ee,\eta,m)\geq m .  \end{equation} }
\end{theorem}

\begin{proof}
The existence of such an index $p(\ee,\eta,m)$ has been proved in Theorem \ref{any_j}.
 For each $m$, let us show that $p(\ee,\eta,m)\geq m$, even if  $[\mathfrak{a},\mathfrak{b}]$ contains points $\eta_0$ such that $\Lambda_m(\eta^0)=\Lambda_{m+1}(\eta^0)$. We divide the proof in several steps depending on $m$.

 \medskip

 {\it  \underline {1.- The case when $m=1$.}} The result $p(\ee,\eta,m)\geq m$ is   self-evident when $m=1$, {and  $[\mathfrak{a},\mathfrak{b}]=[-\pi,\pi]$. In addition, when $\eta=\pm \pi$ the multiplicity  of $\Lambda_1^0(\pm \pi)$ is $2$, and $(\pm \pi, \pi^2)$ are the only points of the limit dispersion curve $\Lambda^0_1(\eta) $ where the multiplicity is greater than $1$.  That is to say,  for
    $\eta=\pm \pi$, the first  eigenvalue of the limit problem is double: $\Lambda_1^0(\pm\pi)=\Lambda_2^0(\pm\pi)=\pi^2$.  The corresponding eigenfunctions  are $U^0_1(-\pi) =e^{-i  \pi  x_1}$ and $U^0_2(-\pi) =e^{i  \pi  x_1}$ when $\eta =-\pi$, while they are  $U^0_1(\pi) =e^{i  \pi  x_1}$ and $U^0_2(\pi) =e^{-i  \pi  x_1}$ when $\eta =\pi$, and we can prove that   $p(\ee,\pm \pi, 2)\geq 2$. In fact,  let us show that there are at least two eigenvalues of \eqref{(11)}--\eqref{(14)} satisfying \eqref{(para1)}. Since both points can be treated in the same way, let us  proceed with $\eta =-\pi$.
  The constructed quasimodes $U_{as,1}^\ee(-\pi)$ and $U_{as,2}^\ee(-\pi)$ satisfy the condition \eqref{(un64as)}. It is self-evident that actually,  $U_{as,1}^\ee(-\pi):=
  U_{as,\pm,0}^\varepsilon( -\pi)$ and $ U_{as,2}^\ee(-\pi)=
 U_{as,+,1}^\varepsilon(-\pi)$, cf.    \eqref{kp0}, \eqref{kpj}, \eqref{(un3)} and \eqref{(un2)}, and, further specifying,
  $$U_{as,1}^\ee(-\pi)=
 \|U_{\pm,0}^\varepsilon(-\pi);{\cal H}^\varepsilon(-\pi)\|^{-1} U^\ee_{\pm,0}(-\pi),\quad    U_{as,2}^\ee(-\pi)=
 \|U_{+,1}^\varepsilon(-\pi);{\cal H}^\varepsilon(-\pi)\|^{-1} U^\ee_{+,1}(-\pi).
 $$

  For any $c_r\geq c$ (cf. \eqref{(un70r)} and \eqref{(un10r)}) let us consider all the eigenvalues $$\{M^\ee_{J}(-\pi) , \cdots, M^\ee_{J+K-1}( -\pi)\} $$ in the interval $[(1+\Lambda^0_1(-\pi))^{-1}-c_r\ee, \, (1+\Lambda^0_1(-\pi))^{-1}+c_r\ee]$  and the corresponding eigenfunctions $\{U^\ee_{J} (\cdot;-\pi) , \cdots, U^\ee_{J+K-1}(\cdot;-\pi) \} $. Using  the  bound       \eqref{(J5e)} in Lemma \ref{Lemma_Visik}, we get
$$
\| U^\varepsilon_{as,l}(-\pi)-\sum_{i=J}^{J+K-1} \alpha_{i,l}^\ee U_i^\varepsilon (\cdot;-\pi)
 ;{\cal H}^\varepsilon(-\pi)\|\leq 2\frac{c}{c_r}, \quad \sum_{i=J}^{J+K-1}\vert\alpha_{i,l}^\ee\vert^2 =1, \quad l=1,2.
$$
Let us show that $\widetilde U_1^\ee := \sum_{i=J}^{J+K-1} \alpha_{i,1}^\ee U_i^\varepsilon (\cdot;-\pi)  $ and   $\widetilde U_2^\ee :=\sum_{i=J}^{J+K-1} \alpha_{i,2}^\ee U_i^\varepsilon (\cdot;-\pi)    $ are linearly independent functions, and consequently, among the sequence $\{M^\ee_{J}( -\pi) , \cdots, M^\ee_{J+K-1}( -\pi)\} $  there are at least   two  eigenvalues of \eqref{(11)}--\eqref{(14)}  with  total multiplicity greater than or equal to $2$.

Indeed,   using the above estimate and \eqref{(un64as)},   we write
\begin{eqnarray}\nonumber
 \left\vert \langle  \widetilde U_{1}^\varepsilon\,,\widetilde U_{2}^\varepsilon \rangle_{-\pi\ee}\right\vert& \leq &
  \left\vert \langle \widetilde U_{1}^\varepsilon\,,\,\widetilde U_{2}^\varepsilon-U_{as,2}^\ee (-\pi) \rangle_{-\pi\ee}\right \vert +
   \left\vert \langle \widetilde U_{1}^\varepsilon  - U_{as,1}^\ee(-\pi) \,,\, U_{as,2}^\ee (-\pi) \rangle_{-\pi\ee}\right \vert  \\ &&   + \left\vert \langle  U_{as,1}^\ee (-\pi) \,,\, U_{as,2}^\ee (-\pi)  \rangle_{-\pi\ee}  \right\vert \leq
 2\frac{c}{c_r} +  2\frac{c}{c_r} + \widetilde C\sqrt{\ee}.  \label{npi}
\end{eqnarray}
 Now, assuming that  $\widetilde U_{1}^\varepsilon$ and $\widetilde U_{2}^\varepsilon$ are two linearly  dependent functions, without any restriction, we can write  $\alpha \widetilde U_1^\ee= \widetilde U_2^\ee $ for some $\alpha \not =0$, $\vert \alpha\vert=1$, and taking the scalar product with $\widetilde U_{1}^\varepsilon$ we have
 $$\vert \alpha \vert
\left\vert \langle  \widetilde U_{1}^\varepsilon\,,\widetilde U_{1}^\varepsilon \rangle_{-\pi\ee}\right \vert  =  1 =
\left\vert \langle  \widetilde U_{2}^\varepsilon\,,\widetilde U_{1}^\varepsilon \rangle_{-\pi\ee}\right \vert    \leq 2\frac{c}{c_r} +  2\frac{c}{c_r} + \widetilde C\sqrt{\ee} .$$
Consequently, it suffices to take $c_r :=c_{r_0}$ and $\varepsilon_0:=\ee_{0,r_0}$ such that  $$4\frac{c}{c_{r_0}} +  \widetilde C\sqrt{ \ee_{0,r_0}} <1$$ to get a contradiction.  Thus, we have also shown that $p(\ee,-\pi,2)\geq 2$ and  $p(\ee,\pi,2)\geq 2$, and we can fix
constants $C_0:=C_{r_0}$ and $\ee_0:=\ee_{0,r_0}$ in \eqref{conm}  to obtain $p(\ee,\eta,1)\geq 1$ $\forall \eta\in [\mathfrak{a},\mathfrak{b}]\equiv[-\pi,\pi]$ while at the points $\eta=\pm \pi$ where    the multiplicity of $\Lambda_1^0(\eta)$ is two, also the index $p(\ee,\pm \pi, 2)$ is  greater than or equal to  $2$. }

\medskip

 {\it  \underline {2.- The case when $m=2$.}}  Let us proceed with $m=2$, in such a way that also $[\mathfrak{a},\mathfrak{b}]=[-\pi,\pi]$, namely,   $H\in(0,1/2)$, cf. Figures \ref{fig_new}  and \ref{fig2} a).  Now, since $\Lambda^0_2(-\pi)= \Lambda^0_1(-\pi)=\pi^2$ and  $\Lambda^0_2(+\pi)= \Lambda^0_1(+\pi)=\pi^2$, we can avoid the two points since the result is the same as for $\Lambda_1^0(\pm \pi)$ and   $p(\ee,\pm \pi, 2)\geq 2$.

For each   $\ee\leq\ee_0$, let us consider the abscises   $\tilde a_{\pm \pi,\ee}$  of the two points in which the  $\ee$-neighborhood of $(1+\Lambda_2^0(\eta))^{-1}$ and the $\ee$-neighborhood of $(1+\Lambda_1^0(\eta))^{-1}$ cut across each other. Further specifying, we consider  the cut points of the lines
\begin{equation}\label{corte}     (1+\Lambda^0_1(\eta))^{-1}-c_{r_0} \ee  \quad \mbox{ and } \quad  (1+\Lambda^0_2(\eta))^{-1}+ c_{r_0}\ee,
 \end{equation}
which we denote by $\tilde a_{-\pi,\ee}$ and  $\tilde a_{\pi,\ee}$,  respectively (see Figures  \ref{inversos} and  \ref{inversoszoom}).

\begin{figure}[ht]
\begin{center}
\resizebox{!}{7cm} {\includegraphics{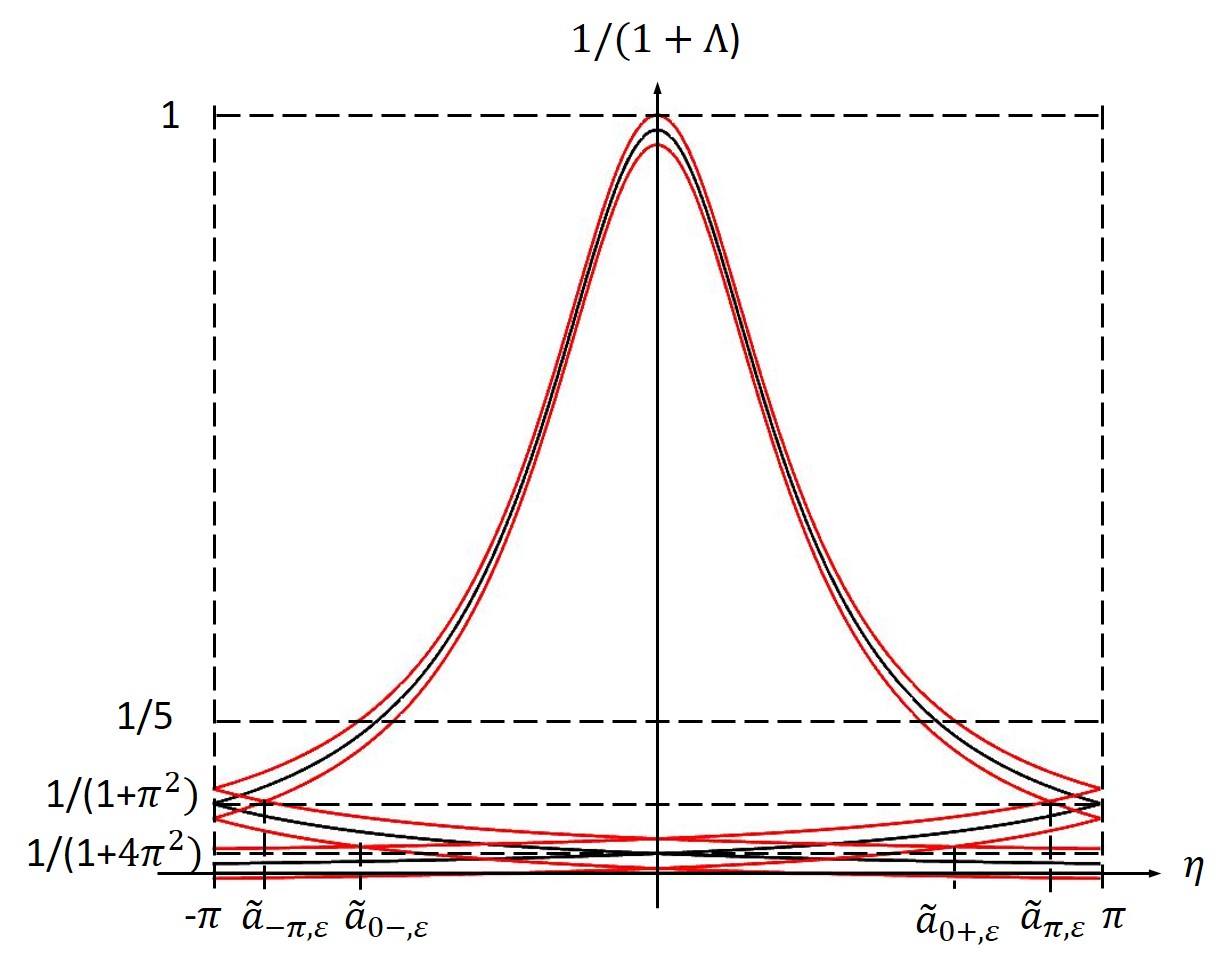}}
\caption{Graphics of the  curves $(1+\Lambda^0_1(\eta))^{-1} $ and $(1+\Lambda^0_2(\eta))^{-1}$  when  $H\in (0,1/2)$.   Associated $\ee$-neighborhoods are   the regions between the   surrounding red lines. }
\label{inversos}
\end{center}
\end{figure}

\begin{figure}[ht]
\begin{center}
\resizebox{!}{6cm} {\includegraphics{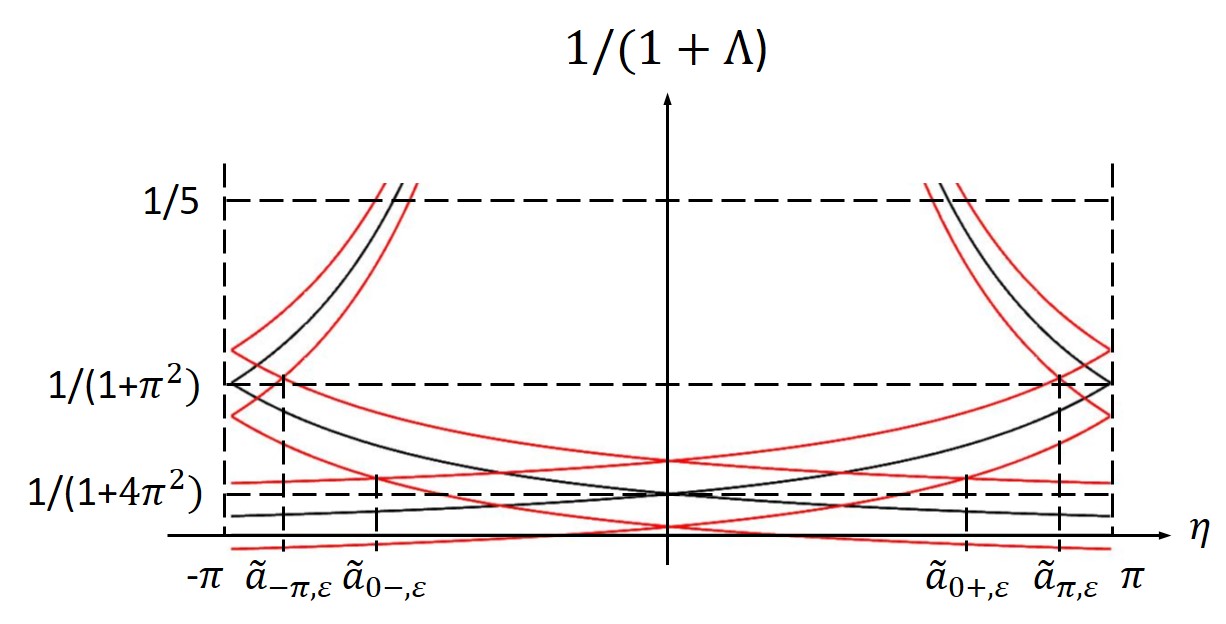}}
\caption{Figure \ref{inversos} after stretching the ordinate axis: the curves $(1+\Lambda^0_1(\eta))^{-1} $ and $(1+\Lambda^0_2(\eta))^{-1}$ and the associated $\ee$-neighborhood between the red lines. }
\label{inversoszoom}
\end{center}
\end{figure}

Consider $\eta$ in the interval $(\tilde a_{-\pi,\ee},\tilde a_{\pi,\ee})$ since the $\ee$-neighborhoods do not intersect,  formulas \eqref{(un70r)}-\eqref{(un10r)} when $j=1$  and  \eqref{(para1)}-\eqref{(un10bis)} provide  two  eigenvalues of \eqref{(11)}--\eqref{(14)} which cannot coincide and the result holds true, namely $p(\ee,\eta,2)\geq 2$  $\forall \eta \in  (\tilde a_{-\pi,\ee},\tilde a_{\pi,\ee})$.

Consider now $\eta$ in the interval $(-\pi, \tilde a_{-\pi, \ee}]$ and $\ee\le \ee_0$.
In this case, both $\ee$-neighborhoods intersect and we take a  larger interval $[(1+\Lambda^0_2(\eta))^{-1}-3c_{ r_0 }\ee, (1+\Lambda^0_2(\eta))^{-1}+3c_{ r_0 }\ee]$ which contains   $[(1+\Lambda^0_1(\eta))^{-1}- c_{ r_0 }\ee, (1+\Lambda^0_1(\eta))^{-1}+ c_{ r_0, }\ee]$ and therefore, it contains at least one  eigenvalue of the operator ${\cal B}^\ee(\eta)$ in \eqref{(J3)}. Let us show that it contains another different eigenvalue of this  operator   and  the result of the theorem is also true for $[\mathfrak{a},\mathfrak{b}]=[-\pi,\pi]$ when  $H\in(0,1/2)$. Here, and below throughout the proof, by different we understand that the total multiplicity of these two eigenvalues is greater than or equal to 2, or equivalently, that they have two corresponding entries in the sequence \eqref{(17)} of eigenvalues of \eqref{(11)}-\eqref{(14)}.

Indeed, formulas \eqref{(un70r)}  and \eqref{(para1)}  provide  two eigenvalues of \eqref{(11)}--\eqref{(14)} which, by now, could coincide. But,  we have constructed two quasimodes $U_{as,1}^\ee(\eta) $ and $U_{as,2}^\ee(\eta)$ satisfying \eqref{(un64as)}.  Let us show that they approach two linear combination of eigenfunctions of  \eqref{(11)}--\eqref{(14)} which are linearly independent functions.
To do this, we consider all the possible eigenvalues of  the operator ${\cal B}^\ee(\eta)$ in the interval $[(1+\Lambda^0_2(\eta))^{-1}- 3c_{ r_0 }\ee, (1+\Lambda^0_2(\eta))^{-1}+ 3c_{ r_0 }\ee]$,  $\{  M^\ee_{i}(\eta)\}_{i=J_\eta}  ^{J_\eta+K_\eta-1}$ and the corresponding eigenfunctions,
and use   \eqref{(J5e)} in Lemma \ref{Lemma_Visik}; we get
$$
\Big\Vert U^\varepsilon_{as,l}(\eta)-\sum_{i=J_\eta}^{J_\eta+K_\eta-1} \alpha_{i,l,\eta}^\ee U_i^\varepsilon (\cdot; \eta)
 ;{\cal H}^\varepsilon(\eta)\Big\Vert \leq 2\frac{c}{ 3c_{r_0}}, \quad \sum_{i=J_\eta}^{J_\eta+K_\eta-1}\vert\alpha_{i,l,\eta}^\ee\vert^2 =1, \quad l=1,2,
$$
and proceed  as above for $\eta= -\pi$,   cf. \eqref{npi}, by contradiction, to obtain the result: namely, we claim that there are at least two different eigenvalues  of the operator ${\cal B}^\ee(\eta)$   in  $[(1+\Lambda^0_2(\eta))^{-1}- 3c_{ r_0 }\ee, (1+\Lambda^0_2(\eta))^{-1}+ 3c_{ r_0 }\ee]$. Then,  we fix
constants $C_0$ and $ \ee_0$ from   $3c_{ r_0 }$  and $\ee_0:=\ee_{r_0}$ in \eqref{conm}, cf. \eqref{(un8)}-\eqref{(un10)},  to obtain $p(\ee,\eta,2)$ greater than or equal to 2, $\forall \eta \in [-\pi,\pi]$.

In addition, when $0<H < 1/2$ (cf. Figures \ref{fig_new} and \ref{fig2}), we can show that at the point  $\eta\equiv 0$,   where  $\Lambda_2^0(0)=\Lambda_3^0(0)=4\pi^2$,  the index  $p(\ee,0,3)$  satisfies  $p(\ee,0,3)\geq 3$.  Indeed, consider  all the  eigenvalues of the operator ${\cal B}^\ee(\eta)$ when $\eta\equiv 0$  in  $[(1+\Lambda^0_2(0))^{-1}- 3c_{ r_0 }\ee, (1+\Lambda^0_2(0))^{-1}+ 3c_{ r_0 }\ee]$.
Formulas \eqref{(un7)}   provide  two eigenvalues of \eqref{(11)}--\eqref{(14)} which, by now, could coincide. But,  we have constructed two quasimodes $U_{as,2}^\ee(0) $ and $U_{as,3}^\ee(0)$ satisfying \eqref{(un64as)}.  Let us show that the two mentioned eigenvalues  are different.

To do this, we consider all the possible eigenvalues    $$\{M^\ee_{J_0}(0) , \cdots, M^\ee_{J_0+K_0-1}(0)\} $$ in the interval $[(1+\Lambda^0_2(0))^{-1}- 3c_{ r_0 }\ee, (1+\Lambda^0_2(0))^{-1}+ 3c_{ r_0 }\ee]$\,     with the  corresponding eigenfunctions $\{U^\ee_{J_0} (\cdot;0), \cdots, U^\ee_{J_0+K_0-1}(\cdot;0)\} $,
and  use the bounds for discrepancies in   Lemma \ref{Lemma_Visik},    cf.  \eqref{(J5e)}; we get
$$
\| U^\varepsilon_{as,l}( 0)-\sum_{i=J_0}^{J_0+K_0-1} \alpha_{i,l,0}^\ee U_i^\varepsilon (\cdot; 0)
 ;{\cal H}^\varepsilon(0)\|\leq 2\frac{c}{ 3c_{r_0}}, \quad \sum_{i=J_0}^{J_0+K_0-1}\vert\alpha_{i,l,0}^\ee\vert^2 =1, \quad l=2,3.
$$
and proceed as  for $\eta= -\pi$, by contradiction, cf. \eqref{npi}, finding  two linearly independent functions
$$\Big\{\sum_{j=J_0}^{J_0+K_0-1} \alpha_{j,2,0}^\ee U_j^\varepsilon ( \cdot; 0)\,, \sum_{j=J_0}^{J_0+K_0-1} \alpha_{j,3,0}^\ee U_j^\varepsilon ( \cdot; 0)\Big\}.$$    Thus, for $\eta\equiv 0$ there are at least two different eigenvalues   of the operator ${\cal B}^\ee(\eta)$, cf. \eqref{(J3)},
  in  the interval $[(1+\Lambda^0_2(0))^{-1}- 3c_{ r_0 }\ee, (1+\Lambda^0_2(0))^{-1}+ 3c_{ r_0 }\ee]$.

Consequently, we obtain $p(\ee,0,2)\geq 2$ and $p(\ee,0,3)\geq 3$. This comes from the fact that the value $\Lambda^\ee_{p(\ee,0,1)}(0)$ cannot coincide with  $\Lambda^\ee_{p(\ee,0,2)}(0)$ or $\Lambda^\ee_{p(\ee,0,3)}(0)$.

\medskip
Now, in the case where $H\geq 1/2$, cf. Figures \ref{fig2}  b) and \ref{fig_new}, considering any interval  $[\mathfrak{a},\mathfrak{b}]=[-\pi,\mathfrak{b}_0]$ or $[\mathfrak{a},\mathfrak{b}]=[ \mathfrak{a}_0, \pi]$ which does not contain the abscises $\eta_0$ of the intersecting points of the dispersion curves $(\eta\pm 2\pi)^2$ and $\eta^2+\pi^2/H^2$, the result of the theorem holds true. Namely $p(\ee,\eta,2)\geq 2$ for any $\eta  \in [-\pi,\mathfrak{b}_0]  \cup  [\mathfrak{a}_0, \pi]$, is a consequence of the results above for $m=1$, $m=2$ and Theorem \ref{any_j}.

\medskip

 {\it  \underline{3.- The case when $m=3$.}}
 When $H<1/2$, the process can be continued for $m=3$, cf. Figures \ref{fig2} a) and  \ref{fig_new}. Let us show that $p(\ee,\eta,3)\geq 3 $ when $\eta\in [\mathfrak{a},\mathfrak{b}]$ with  $0\in[\mathfrak{a},\mathfrak{b}] \subset [ a_0,b_0]\subset [-\pi,\pi]$, $(a_0,b_0)$  being the largest interval which does not contain abscises of intersecting points of the curves  $((\eta+2\pi)^2+1)^ {-1}$ and $((\eta^2+\pi^2/H^2)+1)^{-1}$ nor of $((\eta-2\pi)^2+1)^ {-1}$ and $((\eta^2+\pi^2/H^2)+1)^{-1}$.

Indeed, let us consider the abscises   $\tilde a_{0\pm,\ee}$ of the cut points of the $\ee$-neighborhood of the curves $( \Lambda_2^0(\eta)+1)^{-1}$ and $(\Lambda_3^0(\eta)+1)^{-1}$  near $\eta=0$,  namely, near  the double eigenvalue $(1+\Lambda_2^0(0))^{-1}= (1+\Lambda_3^0(0))^{-1}=1/(1+4\pi^2)$ of ${\cal  B}^\ee(0)$.
 Further specifying, $\tilde a_{0-,\ee}$ and  $\tilde a_{0+,\ee}$   denote   abscisas of  the intersecting points   of the curves
\begin{equation}\nonumber    (1+\Lambda^0_2(\eta))^{-1}-c_{r_0} \ee  \quad \mbox{ and } \quad  (1+\Lambda^0_3(\eta))^{-1}+ c_{r_0} \ee   ,
 \end{equation}
 with $\tilde a_{0-,\ee}<\tilde a_{0+,\ee}$; see Figures \ref{inversos} and \ref{inversoszoom}.

In the case where $\mathfrak{a}>-\pi$ (similarly, $\mathfrak{b}<\pi$),  it is clear that for $\eta \in   [\mathfrak{a}, \tilde a_{0-,\ee})$ (similarly, $\eta \in (\tilde a_{0+,\ee}, \mathfrak{b}]  $)  the   eigenvalues $\Lambda_i^0(\eta)$ with $i=1,2,3$ are simple, \eqref{(un70)} and \eqref{(para1)} guarantee  at least one eigenvalue $(\Lambda_p^\ee(\eta)+1)^{-1}$ in the $\ee$-neighborhood of each $( \Lambda_i^0(\eta)+1)^{-1}$, and since these neighborhoods do not intersect, there are at least three different eigenvalues $\Lambda_{p(\ee,\eta,i)}^\ee(\eta)$, $i=1,2,3$ (their total multiplicity is greater than or equal to 3). Therefore, $p(\ee,\eta,3)\geq 3$ for $\eta \in   [\mathfrak{a}, \tilde a_{0-,\ee})\cup (\tilde a_{0+,\ee}, \mathfrak{b}]  $.  Also,  if $\mathfrak{a}=-\pi$ (similarly, $\mathfrak{b}=\pi$), it has been proved above that in an $\ee$-neighborhood of  $( \Lambda_1^0(\eta)+1)^{-1}$ there are at least two eigenvalues of the operator ${\cal B}^\ee(\eta)$, cf. \eqref{(J3)}, with total multiplicity greater than or equal to 2, and again $p(\ee,\eta,3)\geq 3$.

Next, let us prove that the  result holds for $\eta\in [ \tilde a_{0-,\ee} ,  \tilde a_{0+,\ee}]$.

 As for $\eta=0$, the result $p(\ee,0,3)\geq 3$ has been proved above, in the previous step.  Because of the symmetry,  it suffices to analyze in further detail the value of  $p(\ee,\eta,3)$ for the case where $\eta \in [\tilde a_{0-,\ee}, 0)$.   We use the idea in the previous step for $\Lambda_2^0(\eta)$ with $\eta$ near $-\pi$. That is,
  now, both $\ee$-neighborhoods intersect and we take a  larger interval $[(1+\Lambda^0_3(\eta))^{-1}-3c_{ r_0 }\ee, (1+\Lambda^0_3(\eta))^{-1}+3c_{ r_0 }\ee]$ which contains   $[(1+\Lambda^0_2(\eta))^{-1}- c_{ r_0 }\ee, (1+\Lambda^0_2(\eta))^{-1}+ c_{ r_0, }\ee]$ and therefore, it contains at least one  eigenvalue of the operator ${\cal B}^\ee(\eta)$, cf. \eqref{(J3)}.  Let us show that it contains another different eigenvalue of  this operator.

Indeed, formulas \eqref{(un70)}, with $j=1$,  provide  two eigenvalues \eqref{(11)}--\eqref{(14)} which, by now, could coincide. But,  we have constructed two quasimodes $U_{as,2}^\ee(\eta) $ and $U_{as,3}^\ee(\eta)$ satisfying \eqref{(un64as)}.  Let us derive that they approach to two linear combination of eigenfunctions of  \eqref{(11)}--\eqref{(14)} which are linearly independent functions.
To do this, we consider all the possible eigenvalues of ${\cal B}^\ee(\eta)$ in $[(1+\Lambda^0_3(\eta))^{-1}- 3c_{ r_0 }\ee, (1+\Lambda^0_3(\eta))^{-1}+ 3c_{ r_0 }\ee]$, $\{  M^\ee_{i}(\eta)\}_{i=J'_\eta}^{J'_\eta+K'_\eta-1}$ and the corresponding eigenfunctions,
and use    Lemma \ref{Lemma_Visik},    cf.  \eqref{(J5e)}; we get
$$
\Big\Vert U^\varepsilon_{as,l}(\eta)-\sum_{ i=J'_\eta}^{J'_\eta+K'_\eta-1} \dot{\alpha}_{i,l,\eta}^\ee U_i^\varepsilon (\cdot; \eta)
 ;{\cal H}^\varepsilon(\eta)\Big\Vert\leq 2\frac{c}{ 3c_{r_0}}, \quad \sum_{i=J_\eta}^{J_\eta+K_\eta-1}\vert\dot{\alpha}_{i,l,\eta}^\ee\vert^2 =1, \quad l=2,3.
$$
and proceed as above for $\eta= -\pi$ (cf. \eqref{npi}), by contradiction,   to obtain the result. Thus,  there are at least two different eigenvalues  of  the operator ${\cal B}^\ee(\eta)$  in  $[(1+\Lambda^0_3(\eta))^{-1}- 3c_{ r_0 }\ee, (1+\Lambda^0_3(\eta))^{-1}+ 3c_{ r_0 }\ee]$.  Since, for   $\eta\in [ \tilde a_{0-,\ee} ,  \tilde a_{0+,\ee}]$,  there is another eigenvalue in the $\ee$-neighborhood for  $(\Lambda_1^0(\eta)+1)^{-1}$ which does not cut  those $\ee$-neighborhoods under consideration for  $M_{\pm,1}^0(\eta)$ ,  it follows that we have obtained at least three   eigenvalues of \eqref{(11)}--\eqref{(14)}  with total multiplicity greater than or equal to 3. Again, we fix the
constants $C_0$ from $3c_{r_0}$  and $\ee_0:=\ee_{r_0}$ in \eqref{conm}  to obtain the value $p(\ee,\eta,3)\geq 3$   in the statement of the theorem.

As a matter of fact, $[\mathfrak{a},\mathfrak{b}]\equiv [-\pi,\pi]$ for $H<1/\sqrt{8}$, see Figure \ref{fig_new}.
Therefore, \eqref{conm12} and all the statements of the theorem hold true.
\end{proof}

As a consequence of the proof of Theorem \ref{th_epsetabis} when dealing with small neighborhoods of eigenvalues of the homogenized problem   of multiplicity $2$,  we  can state the following result.

\begin{corollary}\label{cor_doble}
{\it Under the hypothesis of Theorem \ref{th_epsetabis}, assume that  for a certain $m\in \{1,2,3\}$      and a certain $\eta^0  \in[\mathfrak{a},\mathfrak{b}]  $,   $\Lambda_m^0(\eta^0)=\Lambda_{m+1}^0(\eta^0)$. Then, for $\ee\leq \ee_0$, there are  $  a_{\eta^0_-, \ee} $ and $  a_{\eta^0_+,\ee}$ such that
  for $\eta\in [ a_{\eta^0_-, \ee}, a_{\eta^0_+,\ee}]\subset [\mathfrak{a},\mathfrak{b}]$,      at least      two eigenvalues $\Lambda^\ee_{p(\ee,\eta,m)}(\eta) $ and $\Lambda^\ee_{p(\ee,\eta,m )+1}(\eta)$   satisfy  \eqref{conm}, and consequently,
\begin{equation}\label{78c}
\left| \Lambda^{\ee}_{p(\ee,\eta,m)}(\eta)- \Lambda^{0}_{j}(\eta) \right| \leq  C_0\ee ,\quad
\left| \Lambda^{\ee}_{p(\ee,\eta,m)+1}(\eta)- \Lambda^{0}_{j}(\eta) \right| \leq  C_0\ee ,\, \quad \forall \eta \in [ a_{\eta^0_-, \ee}, a_{\eta^0_+,\ee}],\,\, \ee\leq \ee_0,
\end{equation}
while $ j=m, m+1$.  The  edges
  $  a_{\eta^0_-, \ee}$  and $a_{\eta^0_+,\ee} $ of the small interval can be determined from an $\ee$-neighborhood of   $(1+\Lambda_m^0(\eta ))^{-1}$ and $(1+\Lambda_{m+1}^0(\eta ))^{-1}$   and both converge towards $\eta^0$, as $\ee\to 0$. }
 \end{corollary}

 }
\begin{remark}\label{rmark}{\rm Comparing results in Theorem \ref{any_j} and  Theorem \ref{th_epsetabis}, it proves useful to observe that  at the points $\eta_0$ such that $\Lambda_m^0(\eta_0)=\Lambda_{m+1}^0(\eta_0)$, the result in Theorem \ref{any_j} provides at least one eigenvalue
$                                                                                                 \Lambda^{\ee}_{p(\ee,\eta_0,m)}(\eta)$  with $ p(\ee,\eta_0,m)\geq m $ but it does not take into account the multiplicity. That is to say,   it does not ensure that there is another different index $p(\ee,\eta_0,m+1)$   such that $ p(\ee,\eta_0,m+1)\geq m+1 $,   and the same can happen  in a
 small neighborhood of   $\eta_0$. This is provided by Theorem \ref{th_epsetabis}.   The main difficulty in the proof of Theorem \ref{th_epsetabis} arises when we are  in an $\ee$-neighborhood of the collision points of the two dispersion curves
 $(\eta_0, \Lambda_m^0(\eta_0))\equiv(\eta_0, \Lambda_{m+1}^0(\eta^0))$, points in which the multiplicity of each branch $\Lambda_{m}^0(\eta )$   and $\Lambda_{m+1}^0(\eta )$ changes from one to two (see Figure  \ref{bandas}).  In these neighborhoods also  the proof  in  Theorem \ref{th_epsetabis} ensures that there are two different indexes, namely, at least the eigenvalues $\Lambda^{\ee}_{p(\ee,\eta,m)}(\eta)$ and $\Lambda^{\ee}_{p(\ee,\eta,m)+1}(\eta)$ satisfy \eqref{conm} for $\eta$ ranging in a small neighborhood  of $\eta_0$,  cf. Corollary \ref{cor_doble}.

 Also,
in connection with the last step of the proof of Theorem \ref{th_epsetabis}, it is worth noting that,   again, depending on the value of $H<1/\sqrt{8}$,   the reasoning  can  be continued for other values of $m$ even greater than $3$ (cf. Figure \ref{fig_new}), and the statement of the theorem   holds true for further values of $m$. However, the process relies on  the complex  trusses-nodes structure for the curves $(1+\Lambda_m^0(\eta))^{-1}$  and computations become cumbersome, this being the reason for avoiding them in this paper. $\Box$}
\end{remark}

\setcounter{equation}{0} \setcounter{section}{4}
\renewcommand{\theequation}{\arabic{section}.\arabic{equation}}

\section{Convergence rates for low-frequency dispersion curves.}\label{sec4}
In this section, we prove the main result of the paper. We show that the index $p(\ee,\eta,m)\geq m$ found in the previous section coincides with $m$,  and as a consequence,  we provide bounds for discrepancies between the eigenvalues of the sequences \eqref{(17)} and \eqref{(S2)} which are uniform in both parameters $\ee$ and $\eta$.

\begin{theorem}\label{teoDc}
{\it Let  $m$ be fixed, $m$ ranging in $\{1,2,3  \}$ and such that $\Lambda^{0}_{m}(\eta)$ is an   eigenvalue of  problem \eqref{(29)}--\eqref{(25)} in the sequence \eqref{(S2)} having a corresponding  eigenfunction independent of $x_2$ for $\eta$ in an interval $[\mathfrak{a},\mathfrak{b}] \subseteq[-\pi,\pi]$.
\begin{itemize} \item[ i).]   Assume that the interval    $ [\mathfrak{a},\mathfrak{b}]  $   does not contain $\eta^0$ such that $\Lambda_m^0(\eta^0)=\Lambda_{m+1}^0(\eta^0)$,   then,
we show that there  exist positive constants $\ee_0$ and $C_0$ {independent of $\eta$ and $\ee$} such that  the eigenvalue $\Lambda^{\ee}_{m}(\eta)$ of problem  \eqref{(11)}--\eqref{(14)} in the sequence \eqref{(17)} is the only one  meeting the estimate
\begin{equation}\label{78cm}
|\Lambda^{\ee}_{m}(\eta)- \Lambda^{0}_{m}(\eta) |\leq C_0 \ee, \quad \forall  \ee\leq \ee_0,\,\,  \eta    \in [\mathfrak{a},\mathfrak{b}].
\end{equation}
\item[ii).]  If   for some $\eta_0 \in [\mathfrak{a},\mathfrak{b}]  $,   $m$ satisfies $\Lambda_m^0(\eta^0)=\Lambda_{m+1}^0(\eta^0)<\Lambda_{m+2}^0(\eta^0)$,  then, there are    exactly  two eigenvalues $\Lambda^\ee_{m}(\eta) $ and $\Lambda^\ee_{m+1}(\eta)$ meeting the discrepancy  \eqref{78cm}  for    $\eta\in [ a_{\eta^0_-, \ee}, a_{\eta^0_+,\ee}]\cap [\mathfrak{a},\mathfrak{b}]$, where $[ a_{\eta^0_-, \ee}, a_{\eta^0_+,\ee}]$ is a small neighborhood  of $\eta_0$.  The  edges   $  a_{\eta^0_-, \ee}$  and $a_{\eta^0_+,\ee} $  can be determined from an $\ee$-neighborhood of    the dispersion curves $\Lambda_m^0(\eta )$ and $\Lambda_{m+1}^0(\eta )$   and both converge towards $\eta^0$, as $\ee\to 0$.

  \item[iii).] In particular, for $m=1$ and $H>0$,  for $m=2$ and  $H\in(0,1/2)$  and for $m=3$ and $H\in (0,1/\sqrt{8})$ we have
\begin{equation}\label{conm12m} |\Lambda_m^\varepsilon(\eta)-\Lambda_m^0(\eta)|<C_0 \varepsilon,\qquad \mbox{for }   \varepsilon\leq \varepsilon_0, \,  \eta \in  [-\pi,\pi] .  \end{equation}
\end{itemize} }
\end{theorem}
\begin{proof} Let us proceed  with the proof in three steps where we use the bounds in Theorem \ref{any_j} and \ref{th_epsetabis}.

\noindent
{\it \underline{First step:} The case where $\forall \eta\in [\mathfrak{a},\mathfrak{b}]  $,   $\Lambda_m(\eta)<\Lambda_{m+1}(\eta)$. }
Let us show that $p(\ee,\eta,m)$ arising in \eqref{conm_j} verifies $p(\ee,\eta,m)=  m$ when $\eta \in [\mathfrak{a},\mathfrak{b}]  $, and consequently \eqref{78cm} also holds.

We proceed by contradiction, denying \eqref{78cm}, while, on account of Theorem \ref{any_j}, $p(\ee,\eta,m)\geq m$.
This implies that  there exists $\eta^*\in [\mathfrak{a},\mathfrak{b}] $  such that there is $\varepsilon_{\eta^*}\leq \ee_m$ for which $p(\ee_{\eta^*}, \eta^*,m)\geq {m+1}$ and \eqref{78cm} does not hold. By the hypothesis,   $$\Lambda_{m+1}^0 (\eta^*)>\Lambda_m^0(\eta^*).$$
 First of all, we observe that, for such a $\eta^*$, the number of  $\ee_{\eta^*}$ that we can select above must constitute a finite number, because otherwise, we can take a subsequence of $ \ee_{\eta^*,l} \to 0$ as $l\to \infty$ and then, from \eqref{conm_j}, or equivalently from \eqref{conm} we write
\begin{equation}\label{etafijo} \Lambda^{\ee_{\eta^*,l}  }_{m+1}(\eta^*) \le  \Lambda_{p(\ee_{\eta^*,l},  \eta^*,m)}^{\ee_{\eta^*,l}  } (\eta^*)
  \leq \Lambda^0 _{m}(\eta^*)  + C_0 {\ee_{\eta^*,l}  }. \end{equation}
 Now, on account of \eqref{(27)}, taking limits, as $l\to \infty$, we get a contradiction:
$$   \Lambda^0_{m+1} (\eta^*) \leq  \Lambda^0_{m}(\eta^*)   .$$
 Consequently, for each $\eta^*$ for which \eqref{78cm} does not hold, there is at the most a finite number $\ee_{\eta^*,1}, \ee_{\eta^*,2}, \cdots \ee_{\eta^*,k_{\eta^*}}$ for which \eqref{78cm} does not hold. In addition, if there is only one  such  $\eta^*$ taking $\ee^*_m= min(\ee_m, \ee_{\eta^*,1}, \ee_{\eta^*,2}, \cdots \ee_{\eta^*,k_{\eta^*}})$, then inequality  \eqref{78cm}  holds for $\ee\leq \ee^*_m$. The same happens if there is only a finite number of $\eta^*$ for which  \eqref{78cm} does not hold.

Thus,  denying \eqref{78cm} must imply    that there is  at least one   sequence $\{\eta_r^*\}_{r=1}^\infty$ that converge towards some $\widehat  \eta \in [ \mathfrak{a},\mathfrak{b}]$, as $r\to \infty$, such that  \eqref{78cm} is not satisfied for $ \ee_{\eta_r^*,1}, \ee_{\eta_r^*,2}, \cdots \ee_{\eta_r^*,k_{\eta_r^*}}$, $r=1,2, \cdots$ while  \eqref{conm_j} holds.
Without any restriction we can assume that there is also a  sequence of $\ee_{\eta_r^*}$ converging towards zero as $r\to \infty$,
 Indeed, let us explain this latter assertion in further detail.  For the set   $\mathcal{J} :=\{\eta^*\in [\mathfrak{a},\mathfrak{b}] : \mbox{ \eqref{78cm} is not satisfied }\}\subset [-\pi,\pi]$, we consider the associated set of parameters constructed above: $$\mathcal{E}:=\{ \ee_{\eta^*,1}, \ee_{\eta^*,2}, \cdots \ee_{\eta^*,k_{\eta^*}}\}_{\eta^*\in \mathcal{J}}.$$  Either $\mathcal {E}$ has a  lower bound $\ee^{**}_m>0$ or we can extract a sequence   $\{\ee_{\eta_r^*}\}_{r=1}^\infty$ converging towards zero as $r\to \infty$,  each term $\ee_{\eta_r^*}$  being associated with a certain value  $\eta_r^*\in \mathcal{J}$. In the first case,  \eqref{78cm} holds for $\ee\leq \ee_m^*:= \min(\ee_m^{**},\ee_m)$ and the proof is completed. In the second case, since the sequence $\{ {\eta_r^*}\}_{r=1}^\infty$ is bounded from above and from below,   by subsequences,  we can  construct the above  mentioned sequence, still denoted by $r$,
 $$(\eta^*_r,  \ee_{\eta_r^*})\to (\widehat \eta, 0) \hbox{ as } r\to \infty.$$ Let us show that this last assertion leads us to a contradiction.

 To this end, we notice that, as in \eqref{etafijo}, from \eqref{conm} we can write
 \begin{equation}\label{ip}\Lambda^{\ee_{\eta_r^*},  }_{m+1}(\eta_r^*) \leq \Lambda^{\ee_{\eta_r^*}  }_{p(\ee_{\eta_r^*}, \eta_r^*,m)}(\eta_r^*) \leq    \Lambda^0 _{m}(\eta_r^*)  + C_0{\ee_{\eta_r^* }  } ,
 \end{equation}
 for   the corresponding sequence of eigenvalues.  Taking into account Theorem \ref{corollary_convergence}, cf. \eqref{cDr}, and the continuity  of the application \eqref{(18m)}, we take limits  in \eqref{ip},  as $r\to \infty$,  and get    $  \Lambda^0_{m+1}(\widehat  \eta) \leq \Lambda_m^0(\widehat  \eta)$. Therefore, when the limiting eigenvalue  $\Lambda_m^0(\widehat  \eta)$ is simple, as happens by the hypothesis on the dispersion curves and the interval  $[\mathfrak{a},\mathfrak{b}]$,  we get a contradiction. The contradiction process also  shows that there are no  more eigenvalues $\Lambda_p^\ee(\eta)$ satisfying  \eqref{78cm}, since by Theorem \ref{any_j} (equivalently, by the hypothesis) $p$ should be greater than or equal to $m$.

\begin{figure}
\begin{center}
\resizebox{!}{7cm} {\includegraphics{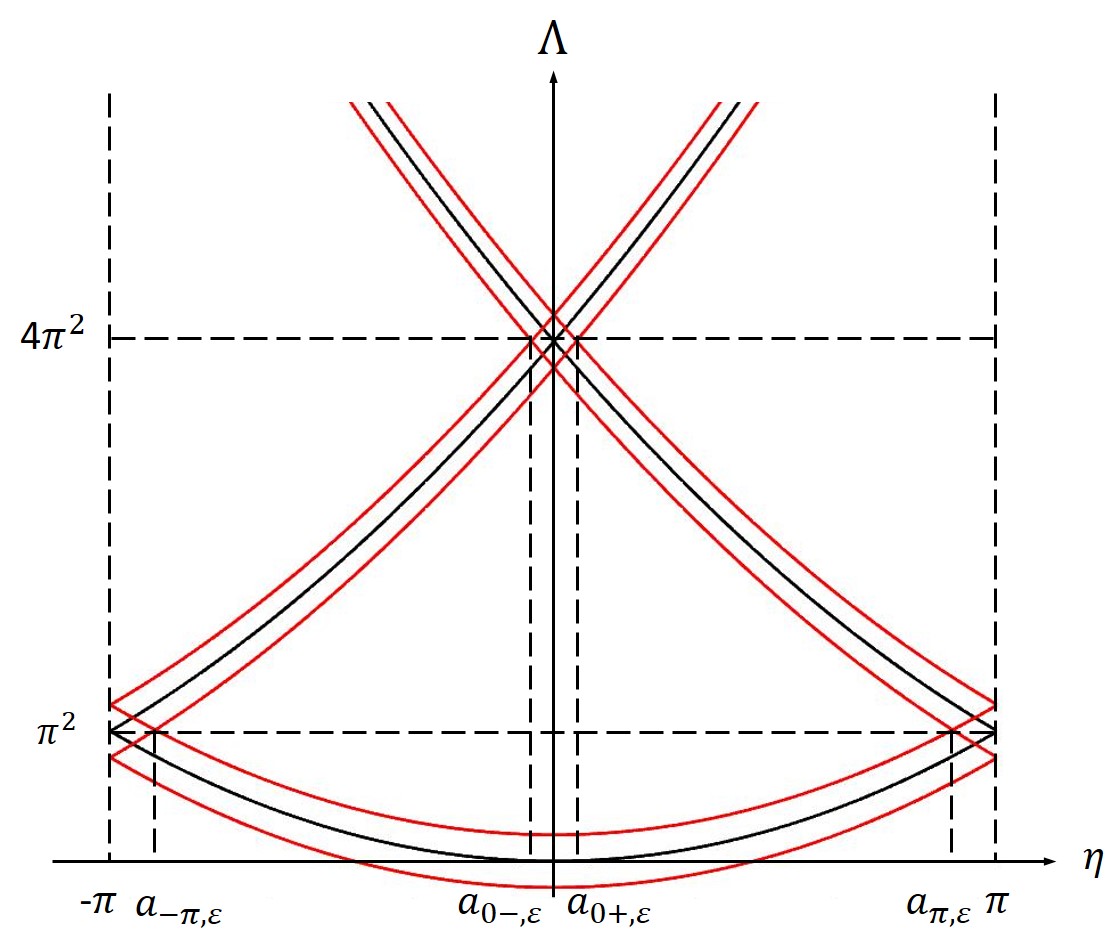}}
\caption{Graphics of the limit dispersion curves  $ \Lambda^0_1(\eta)  $ and $ \Lambda^0_2(\eta) $  when  $H\in (0,1/2)$.   Associated $\ee$-neighborhoods are   the regions between the   surrounding red lines. }
\label{bandas}
\end{center}
\end{figure}

\medskip
\noindent
 {\it \underline{Second step:}  the case where  $\eta^0\in [\mathfrak{a},\mathfrak{b}]$ with $\Lambda_m^0(\eta^0)=\Lambda_{m+1}^0(\eta^0)$.}
Let us start with $m=1$.

By the hypothesis on the limiting dispersion curves, there are only two possibilities for the curves $\Lambda_1^0(\eta)$ and $\Lambda_2^0(\eta)$,   which is to have a  common point for the values of $\eta_0=\pm \pi$, $   [\mathfrak{a},\mathfrak{b}]$  being  either $  [-\pi,\mathfrak{b}]$ when $\eta_0= -\pi$  or $     [\mathfrak{a},\pi]$ when $\eta_0=  \pi$ .  The two nodes are respectively $(-\pi,\pi^2)$ and $(\pi, \pi^2)$ and this holds for any $H>0$.
Let us start by $\eta=-\pi$  and consider the points
$$a_{-\pi,\ee } =-\pi+\alpha_0\ee   \quad \mbox{ and }  \quad a_{ \pi,\ee } =\pi-\alpha_0\ee , \quad \mbox{  with } \alpha_0=\frac{C_{  0}  }{2\pi},\quad \ee\leq \ee_0,$$
which are the abscises of intersecting points of the $\ee$-neighborhood lines of the dispersion curves: $\Lambda_1^0(\eta)+C_0\ee$ and  $\Lambda_2^0(\eta)-C_0\ee$ near $\pm \pi$ (cf. \eqref{corte} to compare).  Since for   $\eta\in [-\pi, a_{-\pi,\ee }]$ the eigenvalues  $\Lambda^{0}_{1}(\eta)$ and  $\Lambda^{0}_{2}(\eta)$ are at a distance $O(\ee)$ between them, it suffices to show
\begin{equation}\label{lp1}
|\Lambda^{\ee}_{1}(\eta)- \Lambda^{0}_{1}(\eta) |\leq C_0 \ee, \quad   |\Lambda^{\ee}_{2}(\eta)- \Lambda^{0}_{1}(\eta) |\leq C_0 \ee, \quad \forall  \ee\leq \ee_0,\,   \eta\in [-\pi, a_{-\pi,\ee }] ,
\end{equation}
with $C_0$ a positive constant independent of $\ee$ and $\eta$.

Considering Theorem \ref{th_epsetabis} (see step 2),    it has been shown that for each $\eta\in [-\pi, a_{-\pi,\ee }]$  there are at least two  eigenvalues  $\Lambda_{p(\ee,\eta,1)}^\ee$  and  $\Lambda_{p(\ee,\eta,1)+1}^\ee$ satisfying $ p(\ee,\eta,1)\geq 1$ and
$$
 |\Lambda_{p(\ee,\eta,1) }^\ee (\eta) - \Lambda^{0}_{1}(\eta) |\leq C_0 \ee, \quad  |\Lambda_{p(\ee,\eta,1)+1}^\ee - \Lambda^{0}_{1}(\eta) |\leq C_0 \ee, \quad \forall  \ee\leq \ee_0,\,   \eta\in [-\pi, a_{-\pi,\ee }] .
$$

Let us assume that  \eqref{lp1}
  does not hold, which means that there is $\ee^*< \ee_0$    for which  we cannot ensure that $p(\ee^* ,\eta_{\ee^*}  ,1)$ is equal to $1$  for some $\eta_{\ee^*}\in[-\pi, a_{-\pi,\ee^* }]$ , namely, $p(\ee^*,\eta_{\ee^*} , 1)\geq 2$  and $p(\ee^*,\eta_{\ee^*} , 1)+1\geq 3$. If   the set of $\ee^*$ for which \eqref{lp1} does not hold is finite, taking the minimum,  we can take the $\ee^{**}$ in such a way  that \eqref{lp1} holds for $\ee\leq\ee^{**}$. Therefore, the worst situation happens when we can extract an infinitesimal  sequence of $\ee_r^*$ and $  \eta_r^*:=\eta_{\ee_r^*}\in  [-\pi, a_{-\pi,\ee_r^* }]$ for which  $p(\ee_r^*,\eta_r^*,  1)\geq 2$  and therefore $p(\ee_r^*,\eta_r^*,  1)+1\geq 3$. If so,
    $(\ee_{\eta_r^*}, \eta^*_r)\to (0, -\pi)$ and, as  in \eqref{ip},   write the corresponding chain of inequalities
  $$\Lambda^{\ee_{\eta_r^*}  }_{3}(\eta_r^*) \leq \Lambda^{\ee_{\eta_r^*}  }_{p(\ee_{\eta_r^*}, \eta_r^*,1)+1}(\eta_r^*) \leq    \Lambda^0 _{1}(\eta_r^*)  + C_0 {\ee_{\eta_r^* }  } $$
 (cf. \eqref{78c}). Taking limits   and applying Theorem \ref{corollary_convergence}, we obtain the following  contradiction
 \begin{equation}\label{23} \Lambda_3^0(- \pi ) \leq   \Lambda^0 _{1}(-\pi )=\Lambda^0 _{2}(-\pi ) < \Lambda_3^0(-\pi).  \end{equation}
Consequently, $p(\ee,\eta, 1)=1$ for $\ee\leq \ee_0,\,   \eta\in [-\pi, a_{-\pi,\ee }] $ and  \eqref{lp1} holds. Since in this interval $\Lambda_1^0(\eta)$ and $\Lambda_2^0(\eta)$ are at distance $O(\ee)$, also  \eqref{78cm} holds for $m=1$ and $m=2$. In addition, there cannot be more eigenvalues $\Lambda_p^\ee(\eta)$  satisfying \eqref{78cm} because of the same argument of contradiction above: indeed,  another eigenvalue should compulsorily  be $\Lambda_3^\ee(\eta)$ for certain $\ee$ and $\eta$ ranging in converging subsequences, which leads to the contradiction \eqref{23}.

In a similar  way, we proceed  with the other possible node of the limiting dispersion curves, that is with the point $(0,4\pi^2)$ when $m=2$ and $H<1/2$. By the hypothesis on the dispersion curves, there is only one possibility for the curves $\Lambda_2^0(\eta)$ and $\Lambda_3^0(\eta)$   to cut across each other. This possibility means to   have a  common point at $\eta_0=0$ and $0\in  [\mathfrak{a},\mathfrak{b}]\subset [-\pi,\pi]$.

We consider the points
 $$a_{0-,\ee}= -\frac{C_0\ee}{ 4\pi}\quad \mbox{ and }\quad a_{0+,\ee}=     \frac{C_0\ee}{ 4\pi},$$
 which are intersecting points of the $\ee$-neighborhood lines of the dispersion curves $\Lambda_2^0(\eta)+C_0\ee$ and  $\Lambda_3^0(\eta)-C_0\ee$ near $\eta=0$, cf. Figure \ref{bandas}.

 Since for $\eta\in [a_{0-,\ee}, a_{0+,\ee }]$
  the eigenvalues  $\Lambda^{0}_{2}(\eta)$ and  $\Lambda^{0}_{3}(\eta)$ are at a distance $O(\ee)$ between them, it suffices to show
\begin{equation}\label{lp2}
|\Lambda^{\ee}_{2}(\eta)- \Lambda^{0}_{2}(\eta) |\leq C_0 \ee, \quad   |\Lambda^{\ee}_{3}(\eta)- \Lambda^{0}_{2}(\eta) |\leq C_0 \ee, \quad \forall  \ee\leq \ee_0,\,    \eta\in [a_{0-,\ee}, a_{0+,\ee }] .
\end{equation}

Also, in Theorem \ref{th_epsetabis} (see step 3, and \eqref{78c}), it has been shown that for each $\eta\in [a_{0-,\ee}, a_{ 0+,\ee }] $  there are at least two   eigenvalues  $\Lambda_{p(\ee,\eta,2)}^\ee(\eta)$  and  $\Lambda_{p(\ee,\eta,2)+1}^\ee(\eta)$ satisfying $ p(\ee,\eta,2)\geq 2$ and
$$
 |\Lambda_{p(\ee,\eta,2) }^\ee (\eta) - \Lambda^{0}_{2}(\eta) |\leq C_0 \ee, \quad   |\Lambda_{p(\ee,\eta,2)+1}^\ee(\eta) - \Lambda^{0}_{2}(\eta) |\leq C_0 \ee, \quad \forall  \ee\leq \ee_0,\,   \eta\in  [a_{0-,\ee}, a_{0+,\ee }] .
$$
  Let us assume that  \eqref{lp2}
  does not hold, which  means that there is $\ee^*$ and $\eta_{\ee^*}\in [a_{0-,\ee^*}, a_{0+,\ee^* }]$ for which  we cannot ensure that $p(\ee^*, \eta_{\ee^*},  2)$ is equal to $2$  , namely, $p(\ee^*, \eta_{\ee^*},  2)\geq 3$  and $p(\ee^*, \eta_{\ee^*},  2)+1\geq 4$. As in the reasoning above (in the first step), the worst situation happens when we can extract a sequence
   $(\ee_{\eta_r^*}, \eta^*_r)\to (0, 0)$ and, as   in \eqref{ip}, we   write
  $$\Lambda^{\ee_{\eta_r^*}  }_{4}(\eta_r^*) \leq \Lambda^{\ee_{\eta_r^*}  }_{p(\ee_{\eta_r^*}, \eta_r^*,2)+1}(\eta_r^*) \leq    \Lambda^0 _{2}(\eta_r^*)  + C_0 {\ee_{\eta_r^* }  }. $$
 Taking limits  and applying  Theorem \ref{corollary_convergence}, we obtain the following  contradiction
\begin{equation}\label{34} \Lambda_4^0(0 ) \leq  \Lambda^0 _{2}(0 ) = \Lambda^0 _{3}(0 ) < \Lambda_4^0(0).  \end{equation}
Consequently, $p(\ee,\eta,  2)=2$ for $\ee\leq \ee_0,\,   \eta\in  [a_{0-,\ee }, a_{0+,\ee  }] $ and  \eqref{lp2} holds as well as \eqref{78cm} for $m=2$ and $m=3$. In addition, there cannot be more eigenvalues $\Lambda_p^\ee(\eta)$  satisfying \eqref{78cm} because of the same argument of contradiction above: indeed,  another eigenvalue should compulsorily  be $\Lambda_4^\ee(\eta)$ for certain $\ee$ and $\eta$ ranging in subsequences, which leads to the contradiction \eqref{34}.

It should be noted that when numbering the eigenvalues \eqref{(31)} in the sequence \eqref{(S2)}, the eigenvalue number $m$ can change depending on the values of $\eta\in [-\pi,\pi]$ and   $H>0$. The rest of the node points when $m\geq 2$ can be treated in a similar way,  and we avoid introducing here  the cumbersome computations.
Also, as a matter of fact, if $H<1/\sqrt{8}$ the process can be continued  for values of $m$ greater than $3$.

\medskip
\noindent
 {\it \underline{Third step}.} Let us show \eqref{conm12m} for the different values of $m$. Now, from the previous steps of this proof, we note that while   $\eta\in [-\pi,0]$,  the  points where estimate \eqref{conm12m} has  not been  proved yet are those for $\eta\in [ a_{-\pi,\ee} , \mathfrak{a}]$ when  $m=1,2$, or $\eta \in [\mathfrak{a},  a_{0-,\ee} ]$ when $m=2,3$, with any $\mathfrak{a} $ such that $-\pi<\mathfrak{a}<0$.  The same can be said while  $\eta\in [0,\pi]$,  with  $\eta\in [\mathfrak{b}, a_{\pi,\ee} ]$ and  $m=1,2$, or $\eta \in [   a_{0+ ,\ee} , \mathfrak{b}]$ when $m=2,3$, with any $\mathfrak{b} $ such that $0<\mathfrak{b}<\pi$.

Let us start with $m=1,2$, and $\eta\in [-\pi,\mathfrak{a}]\supsetneqq [ a_{-\pi,\ee} , \mathfrak{a}]$.   Theorem \ref{th_epsetabis} ensures that  \eqref{conm} holds  for $m=2$ with $p(\ee,\eta,2)\geq 2$. We can apply the same contradiction argument used above  which leads to a  certain sequence
$$(\eta^*_r,  \ee_{\eta_r^*})\to (\widehat \eta, 0) \hbox{ as } r\to \infty, \quad \mbox{ for some }\widehat  \eta \in  [-\pi,\mathfrak{a}],$$  for which  $p( \ee_{\eta_r^*},\eta_r^*,2)\geq 3$. Actually, taking limits in
$$\Lambda^{\ee_{\eta_r^*}  }_{3}(\eta_r^*) \leq \Lambda^{\ee_{\eta_r^*}  }_{p(\ee_{\eta_r^*}, \eta_r^*,1) }(\eta_r^*) \leq    \Lambda^0 _{2}(\eta_r^*)  + C_0 {\ee_{\eta_r^* }  } $$
 cf. \eqref{cDr}, we obtain
 $\Lambda_3^0(\widehat  \eta ) \leq   \Lambda^0 _{2}(\widehat  \eta )$ which is a contradiction for   $\widehat  \eta  \in  [  -\pi  , \mathfrak{a}] \subset [-\pi, 0) .   $ Thus, \eqref{conm12m} holds  for $m=2$  and $\eta \in  [  -\pi  , \mathfrak{a}] \subset [-\pi, 0) .   $

 For $m=1$, we observe that the construction of the endpoints in the second step applies if we replace the abscise $a_{-\pi,\ee} $ by $\widetilde a_{-\pi,\ee}$  arising in the second step of the proof of Theorem \ref{th_epsetabis} (see Figures \ref{inversos}, \ref{inversoszoom} and \ref{bandas}). The construction in Theorem \ref{th_epsetabis} ensures that, for $\eta\in [\widetilde a_{-\pi,\ee}, \mathfrak{a}]$, there are at least two different eigenvalues:  $\Lambda_{p (\ee,\eta,1)}^\ee(\eta)$   satisfying \eqref{conm12}  with $m=1$,  and $\Lambda_{p(\ee,\eta,2)}^\ee (\eta)$ satisfying \eqref{conm12}  with   $m=2$. Since we have proved above that $p(\ee,\eta,2)=2$,  also $p (\ee,\eta,1)=1$.   Thus,  \eqref{conm12m} also holds for $m=1$ and $\eta \in [  -\pi  , \mathfrak{a}]$, and therefore we have obtained  that for $H<1/2$ bounds \eqref{conm12m} hold for $m=1,2$.

 We proceed in a similar way in the interval $\eta \in [\mathfrak{a},  a_{0-,\ee} ]$ when $m=2,3$, arguing by contradiction for a $p( \ee_, \eta,3)\geq 4$ and then,   combining the second step   of the proof of this theorem with the third step in the proof of  Theorem  \ref{th_epsetabis}. Thus, we obtain  that for $H<1/\sqrt{8}$ bounds \eqref{conm12m} hold for $m=1,2,3$.

Therefore, the result of the last statement   holds and the theorem is proved.
\end{proof}

{

\end{document}